\documentclass[thmsa,a4paper,11pt]{article}%
\usepackage{amssymb}
\usepackage{amsmath}
\usepackage{sw20jart}
\usepackage{graphicx}
\usepackage{amsfonts}%
\newtheorem{thm}[theorem]{Th\'eor\`eme}
\newtheorem{defin}{D\'efinition}[section]

\newtheorem{lemme}[theorem]{Lemme}
\newtheorem{exem}{Exemple}[section]

\newtheorem{remarque}{Remarque}
\newtheorem{Convention}{Convention}
\newtheorem{Conventions}[Convention]{Conventions}
\newtheorem{remarques}[remarque]{Remarques}

\newtheorem{question}[exercise]{Question}
\begin{document}

\title{Fonction constante et d\'{e}riv\'{e}e nulle~: un r\'{e}sultat si trivial...}
\author{A.\ Delcroix$^{1,2}$ et C.\ Silvy$^{2}$\\(1) Laboratoire AOC\\Universit\'{e} Antilles-Guyane\\UFR Sciences, Campus de Fouillole\\97159 Pointe \`{a} Pitre Cedex, Guadeloupe\\(2) Centre de Recherches et de Ressources en Education et Formation\\IUFM\ de Guadeloupe\\Morne Ferret -- BP\ 517 \\97178\ Abymes Cedex, Guadeloupe}
\date{03 septembre 2008}
\maketitle

\begin{quotation}
On peut remarquer \`{a} ce sujet que si $f^{\prime}$ n'est pas continue, ce
n'est pas non plus une fonction \emph{tr\`{e}s} sauvage (...)

Roger Godement (Analyse math\'{e}matique 1, Springer 2\`{e} \'{e}dition, 2001)
\end{quotation}

\section{Introduction}

Dans le cadre de l'analyse d'une restitution organis\'{e}e de connaissance
donn\'{e}e au baccalaur\'{e}at (s\'{e}rie S, Antilles-Guyane 2006)
\cite{SilvyDelcroix1}, la r\'{e}flexion sur la caract\'{e}risation des
fonctions constantes d\'{e}finies sur un intervalle par leur d\'{e}riv\'{e}e
conduit \`{a} s'interroger sur la \textquotedblleft
d\'{e}montrabilit\'{e}\textquotedblright\ \'{e}ventuelle de cette
propri\'{e}t\'{e} en classe de terminale scientifique.\ La d\'{e}monstration
de
\[
f\ \text{constante}\Rightarrow f\text{ d\'{e}rivable et }f^{\prime
}\ \text{nulle}%
\]
\'{e}tant consid\'{e}r\'{e}e comme une cons\'{e}quence imm\'{e}diate de la
d\'{e}finition de la d\'{e}riv\'{e}e, c'est sur la r\'{e}ciproque que nous
nous concentrons. La d\'{e}monstration qui vient le plus souvent \`{a}
l'esprit fait appel \`{a} l'in\'{e}galit\'{e} (ou \`{a} l'\'{e}galit\'{e}) des
accroissements finis.\ Elle se situe aux bordures mouvantes du programme des
classes scientifiques. Nous avons essay\'{e} -\thinspace en nous situant \`{a}
dire vrai plus dans le programme des classes pr\'{e}paratoires aux grandes
\'{e}coles o\`{u} d'un L1 sientifique\thinspace- de recenser diff\'{e}rentes
d\'{e}monstrations de cette caract\'{e}risation (not\'{e}e FCD dans la
suite).\ Nous avons choisi d'exclure de notre \'{e}tude les d\'{e}monstrations
faisant appel au calcul int\'{e}gral qui, sous des hypoth\`{e}ses convenables,
permettent une d\'{e}monstration imm\'{e}diate. (On retrouvera cependant des
commentaires sur cette question dans la sous partie \ref{FCD2008SPAnFonc}.)
Dans la pratique nous \'{e}tendons notre analyse au principe de Lagrange,
liant le sens de variation de la fonction au signe de la d\'{e}riv\'{e}e
(caract\'{e}risation not\'{e}e SVD dans la suite).\ En effet, la
caract\'{e}risation FCD peut \^{e}tre facilement vue comme une cons\'{e}quence
de la caract\'{e}risation SVD, qui devient alors un r\'{e}sultat cl\'{e}. Pour
cette derni\`{e}re, et pour une fonction d\'{e}rivable, l'implication
\[
f\ \text{croissante (\emph{resp}. d\'{e}croissante)}\Rightarrow f^{\prime
}\ \text{positive (\emph{resp}. n\'{e}gative)}%
\]
\'{e}tant de nouveau cons\'{e}quence imm\'{e}diate des d\'{e}finitions, nous
nous int\'{e}ressons \'{e}galement essentiellement \`{a} sa
r\'{e}ciproque.\smallskip

\noindent Cependant, notons une diff\'{e}rence de nature entre ces deux
r\'{e}sultats. La caract\'{e}risation SVD est \'{e}troitement reli\'{e}e aux
fonctions de la variable et \`{a} valeurs r\'{e}elles, en fait \`{a} la
relation d'ordre sur $\mathbb{R}$. La caract\'{e}risation FCD s'inscrit dans
une question math\'{e}matique plus vaste, celle du noyau de l'op\'{e}rateur
lin\'{e}aire qui \`{a} une fonction associe sa d\'{e}riv\'{e}e.\ C'est une
question importante sur tout espace muni d'une d\'{e}rivation. (Par exemple,
sur l'espace des distributions de Schwartz, la question de la
caract\'{e}risation des distributions de d\'{e}riv\'{e}e nulle reste un
exercice relativement difficile lorsqu'on aborde la th\'{e}orie, voir la sous
partie \ref{FCD2008SPAnFonc}.)\smallskip

\noindent Par ailleurs, ce parcours autour des caract\'{e}risations FCD et SVD
nous a conduit \`{a} remarquer qu'\`{a} la base de ce r\'{e}sultat figurent
les propri\'{e}t\'{e}s essentielles du corps des nombres
r\'{e}els.\ Rappelons, pour faciliter les r\'{e}f\'{e}rences ult\'{e}rieures,
les 5 propri\'{e}t\'{e}s \'{e}quivalentes suivantes (Voir \cite{Rogalski1}
pour des d\'{e}veloppements sur ce sujet)~:

\noindent$\left(  BS\right)  $ \emph{propri\'{e}t\'{e} de la borne
sup\'{e}rieure}~: toute partie de $\mathbb{R}$ non vide major\'{e}e admet une
borne sup\'{e}rieure~;\newline$\left(  MB\right)  $ \emph{convergence des
suites monotones born\'{e}es}~: toute suite r\'{e}elle monotone born\'{e}e est
convergente~;\newline$(SE)$ \emph{propri\'{e}t\'{e} des segments
embo\^{\i}t\'{e}s}~: une suite d\'{e}croissante (pour l'inclusion)
d'intervalles r\'{e}els ferm\'{e}s born\'{e}s de longueurs tendant vers $0$
poss\`{e}de une intersection non vide, r\'{e}duite \`{a} un point~;\newline%
$(CC)$ \emph{compl\'{e}tude s\'{e}quentielle}~: toute suite r\'{e}elle de
Cauchy est convergente~;\newline$(BW)$ \emph{propri\'{e}t\'{e} de
Bolzano-Weierstrass}~: de toute suite r\'{e}elle born\'{e}e, on peut extraire
une suite convergente.\smallskip

\noindent Comme on le verra dans les parties \ref{FCD2008DemoAF} et
\ref{FCDPartieConnexite}, les d\'{e}monstrations bas\'{e}es sur
l'in\'{e}galit\'{e} des accroissements finis (ou sur des m\'{e}thodes issues
de d\'{e}monstrations de cette in\'{e}galit\'{e}) utilisent plut\^{o}t la
propri\'{e}t\'{e} $(BS)$ et la connexit\'{e} des intervalles r\'{e}els.\ Celle
bas\'{e}e sur l'\'{e}galit\'{e} des accroissements finis se r\'{e}f\`{e}re
indirectement \`{a} la propri\'{e}t\'{e} $(BW$), utilis\'{e}e pour
d\'{e}montrer les propri\'{e}t\'{e}s des images de segments par les fonctions
continues. Enfin, une m\'{e}thode bas\'{e}e sur un principe de dichotomie
(voir partie \ref{FCD2008DemoDicho}) se rattache \`{a} la propri\'{e}t\'{e}
$(SE)$ (puisque le principe de dichotomie d\'{e}coule facilement de
$(SE)$).\ Notons que cette derni\`{e}re m\'{e}thode pourrait \emph{a priori}
faire l'objet d'un probl\`{e}me en classe de terminale S, l'ensemble de ses
ingr\'{e}dients (principe de dichotomie, suites convergentes,...) figurant
dans le programme des classes scientifiques du secondaire.\smallskip

\noindent Au del\`{a} de ce simple recensement de techniques de
d\'{e}monstration, nous avons souhait\'{e} les mettre en rapport (partie
\ref{FCD2008ComMath}) avec des notions math\'{e}matiques sous jacentes
(propri\'{e}t\'{e} des pentes, stricte d\'{e}rivabilit\'{e} auquel un des
lemmes utilis\'{e} dans la d\'{e}monstration par dichotomie fait
immanquablement penser).\ Nous remettons \'{e}galement en place
l'\'{e}quivalence entre le th\'{e}or\`{e}me de Rolle, celui des accroissements
finis et, ce qui nous semble moins utilis\'{e}, l'ensemble form\'{e} par le
th\'{e}or\`{e}me de Darboux et la caract\'{e}risation SVD. Nous apportons
quelques commentaires sur la place de la caract\'{e}risation FCD en analyse,
lorsque la d\'{e}riv\'{e}e devient un op\'{e}rateur lin\'{e}aire et la
caract\'{e}risation FCD la recherche d'un noyau. Fort de ces \'{e}l\'{e}ments,
nous construisons alors le site math\'{e}matique
\cite{DuchetErdogan1,SilvyDelcroix1} de la caract\'{e}risation SVD.

\begin{Conventions}
~\newline$\left(  i\right)  $~Lorsqu'on parle d'un intervalle $[a,b]$ dans la
suite, il est sous entendu que $a$ et $b$ sont deux r\'{e}els tels que
$a<b$.~\newline$\left(  ii\right)  $~Nous avons choisi de nous restreindre
(sauf mention explicite du contraire) \`{a} des fonctions d\'{e}finies sur un
segment $[a,b]$, continues sur ce segment et d\'{e}rivables sur $\left]
a,b\right[  $ pour mettre notre \'{e}tude dans le cadre d'hypoth\`{e}ses
classiques en d\'{e}but de premier cycle d'\'{e}tudes sup\'{e}rieures.
\end{Conventions}

\section{D\'{e}monstrations des caract\'{e}risations FCD et SVD par des
propri\'{e}t\'{e}s d'accrois\-sements finis\label{FCD2008DemoAF}}

\subsection{Trois in\'{e}galit\'{e}s des accroissements finis}

Nous commencons par rappeler, pour le moment sans d\'{e}monstration, trois
variantes de l'in\'{e}galit\'{e} des accroissements finis.

\begin{proposition}
\label{FCD2008LemmeACFClassique}(\textbf{In\'{e}galit\'{e} des accroissements
finis classique, IAF}) Soit $f:[a,b]\rightarrow\mathbb{R}$ continue sur
$[a,b]$ et d\'{e}rivable en tout point de $]a,b[$ et $k$ un r\'{e}el positif
tels que~: $\forall x\in\left]  a,b\right[  ,\ \ $%
$\vert$%
$f^{\prime}(x)|\leq k$. Alors%
\[
|f(b)-f(a)|\leq k(b-a).
\]

\end{proposition}

\noindent Cette premi\`{e}re version est souvent privil\'{e}gi\'{e}e, puisque
c'est la version par excellence g\'{e}n\'{e}ralisable au cas des fonctions
\`{a} valeurs dans un espace norm\'{e}.\ Dans le cycle terminal du secondaire,
la formule suivante a pu \^{e}tre pr\'{e}f\'{e}r\'{e}e.

\begin{proposition}
\label{FCD2008LemmeACFCT}(\textbf{In\'{e}galit\'{e} des accroissements finis,
IAF'}) Soit $f:[a,b]\rightarrow\mathbb{R}$, continue sur $[a,b]$ et
d\'{e}rivable en tout point de $]a,b[$ et $m,\ M$ des r\'{e}els positifs tels
que~: $\forall x\in\left]  a,b\right[  ,\ \ m\leq f^{\prime}(x)\leq M$. Alors%
\[
m(b-a)\leq f(b)-f(a)\leq M(b-a).
\]

\end{proposition}

\noindent Ces deux formules conduisent \`{a} une d\'{e}monstration
imm\'{e}diate de la caract\'{e}risation FCD.\ En effet, sur tout intervalle
$\left[  a,\beta\right]  \subset\lbrack a,b]$, les hypoth\`{e}ses de la
proposition \ref{FCD2008LemmeACFClassique} (resp. \ref{FCD2008LemmeACFCT})
sont satisfaites avec $k=0$ (resp. $m=M=0$), conduisant \`{a} $f(a)=f(\beta)$
pour tout $\beta\in\lbrack a,b]$. Pour la troisi\`{e}me version de
l'in\'{e}galit\'{e}, nous suivons la d\'{e}marche propos\'{e}e par le
\textit{Cours de math\'{e}matiques sp\'{e}ciales }(E.\ Ramis \textit{et alii},
3e \'{e}dition, 1991).

\begin{proposition}
\label{FCD2008LemmeACFRamis}(\textbf{In\'{e}galit\'{e} des accroissements
finis g\'{e}n\'{e}ralis\'{e}e, IAFG}) Soit $f,g:[a,b]\rightarrow\mathbb{R}$,
continues sur $[a,b]$ et d\'{e}rivables \`{a} droite en tout point de $]a,b[$,
telles que~: $\forall x\in\left]  a,b\right[  ,\ \ $%
$\vert$%
$f_{d}^{\prime}(x)|\leq g_{d}^{\prime}(x)$ ($f_{d}^{\prime}$ et $g_{d}%
^{\prime}$ d\'{e}signent les d\'{e}riv\'{e}es \`{a} droite). Alors%
\[
|f(b)-f(a)|\leq g(b)-g(a).
\]

\end{proposition}

\noindent En fait, dans l'ouvrage cit\'{e}, la d\'{e}monstration est faite
pour une fonction $f$ \`{a} valeurs dans un espace norm\'{e}, mais nous avons
choisi de nous restreindre aux fonctions \`{a} valeurs r\'{e}elles\footnote{La
d\'{e}monstration est tr\`{e}s similaire \`{a} celle que nous d\'{e}veloppons
dans la partie \ref{FCDPartieConnexite} pour la caract\'{e}risation FCD par un
argument de connexit\'{e}. Elle repose, pour beaucoup, sur la connexit\'{e}
des intervalles r\'{e}els et sur la propri\'{e}t\'{e} $\left(  BS\right)  $.}.
Ce r\'{e}sultat sert \`{a} \'{e}tablir qu'une fonction $f:[a,b]\rightarrow
\mathbb{R}$, d\'{e}rivable \`{a} droite en tout point de $]a,b[$ est
$k$-lipschitizenne sur $[a,b]$ si, et seulement si,%
\[
\forall x\in\left]  a,b\right[  ,\ \ |f_{d}^{\prime}(x)|\leq k.
\]
On en d\'{e}duit alors la caract\'{e}risation FCD, pour une fonction
$f$\ d\'{e}rivable \`{a} droite, par les \'{e}quivalences%
\[
f\ \text{constante}\Leftrightarrow f\ \text{est\textrm{\ }}%
k\text{-lipschitzienne\ de\ constante}\ k=0\Leftrightarrow f_{d}^{\prime}=0.
\]

\begin{remarques}
\label{FCD2008RemSVD}~\newline$\left(  i\right)  $~Le lemme
\ref{FCD2008LemmeACFRamis} sert aussi \`{a} \'{e}tablir le principe de
Lagrange (caract\'{e}risation SVD), ici pris comme le lien entre le sens de
variation d'une fonction $g:[a,b]\rightarrow\mathbb{R}$ d\'{e}rivable \`{a}
droite sur $\left]  a,b\right[  $ et le signe de sa d\'{e}riv\'{e}e \`{a}
droite. En effet, supposons par exemple $g_{d}^{\prime}(t)$ positive sur
$\left]  a,b\right[  $.\ En faisant $f=0$ dans le lemme
\ref{FCD2008LemmeACFRamis}, on obtient $0\leq g(b)-g(a)$, d'o\`{u} on
d\'{e}duit le r\'{e}sultat.\newline$\left(  ii\right)  $~En supposant
$f_{d}^{\prime}$ major\'{e}e sur $]a,b[$ par $k\geq0$ et en faisant $g\left(
x\right)  =kx$, on retrouve l'in\'{e}galit\'{e} classique des accroissements
finis~:%
\[
|f(b)-f(a)|\leq k(b-a).
\]

\end{remarques}

\subsection{Une propri\'{e}t\'{e} centrale~: la majoration des accroissements}

Pour des raisons d'homog\'{e}n\'{e}it\'{e} de l'expos\'{e} nous supposons dans
ce paragraphe que la fonction $f$ (et le cas \'{e}cheant la fonction $g$) sont
d\'{e}rivables sur $\left]  a,b\right[  $. Commencons par introduire une
propri\'{e}t\'{e} de majoration des accroissements.

\begin{lemme}
(\textbf{Majoration des accroissements, MAJA}) Soit $f:[a,b]\rightarrow
\mathbb{R}$, continue sur $[a,b]$ et d\'{e}rivable en tout point de $]a,b[$ et
$M$ un r\'{e}el positif tels que~: $\forall x\in\left]  a,b\right[
,\ \ f^{\prime}(x)\leq M$. Alors%
\[
f(b)-f(a)\leq M(b-a).
\]

\end{lemme}

\noindent Remarquons que~:\vspace{-0.05in}

\begin{enumerate}
\item[$(i)$] l'in\'{e}galit\'{e} des accroissements finis classiques
entra\^{\i}ne la propri\'{e}t\'{e} IAF' (elle lui est donc \'{e}quivalente, la
r\'{e}ciproque s'obtenant en prenant $m=-k$ et $M=k$),\vspace{-0.05in}

\item[$(ii)$] l'in\'{e}galit\'{e} des accroissements finis
g\'{e}n\'{e}ralis\'{e}s entra\^{\i}ne la propri\'{e}t\'{e} MAJA,\vspace
{-0.05in}

\item[$(iii)$] la propri\'{e}t\'{e} MAJA entra\^{\i}ne la caract\'{e}risation
SVD,\vspace{-0.05in}

\item[$(iv)$] la caract\'{e}risation SVD entra\^{\i}ne l'in\'{e}galit\'{e} des
accroissements finis g\'{e}n\'{e}ralis\'{e}s (IAFG).
\end{enumerate}

\noindent En effet pour le $\left(  i\right)  $, on applique, sous les
hypoth\`{e}ses de la propri\'{e}t\'{e} IAF', la propri\'{e}t\'{e} IAF \`{a}
$f_{1}(x)=f(x)-mx$.\ En remarquant que, pour tout $x\in\left]  a,b\right[  $,
on a $0\leq f_{1}^{\prime}(x)\leq M-m$ (d'o\`{u} $0\leq\left\vert
f_{1}^{\prime}(x)\right\vert \leq M-m$), il vient $\left\vert f_{1}%
(b)-f_{1}(a)\right\vert \leq\left(  M-m\right)  \left(  b-a\right)  $.\ Puis,
en particulier,
\[
f_{1}(b)-f_{1}(a)=f(b)-f(a)-m\left(  b-a\right)  \leq\left(  M-m\right)
\left(  b-a\right)
\]
D'o\`{u} $f(b)-f(a)\leq M\left(  b-a\right)  $. On renouvelle l'op\'{e}ration
avec $f_{2}(x)=Mx-f(x)$ pour obtenir l'autre in\'{e}galit\'{e} de IAF'. Pour
le $\left(  ii\right)  $, on applique la propri\'{e}t\'{e} IAFG au couple
$\left(  0,M-f\right)  $. Le $\left(  iii\right)  $ est imm\'{e}diat~: si, par
exemple, $f^{\prime}$ est positive sur $\left]  a,b\right[  $, on applique la
propri\'{e}t\'{e} MAJA au couple $\left(  -f,0\right)  $ pour obtenir
$f(b)\geq f(a)$ et conclure \`{a} la croissance de $f$. Enfin le $\left(
iv\right)  $ repose sur la r\'{e}\'{e}criture sans valeur absolue de
l'hypoth\`{e}se de majoration de la d\'{e}riv\'{e}e de $f$, id est
\[
\forall t\in\left]  a,b\right[  ,\ \ -g^{\prime}(t)\leq f^{\prime}(t)\leq
g^{\prime}(t).
\]
La caract\'{e}risation SVD entra\^{\i}ne alors que la fonction $g-f$
(resp.\ $f+g$) est croissante.\ D'o\`{u}
\begin{align*}
\left(  g-f\right)  (a)  &  \leq\left(  g-f\right)  (b)\text{ id est
}f(b)-f(a)\leq g(b)-g(a)\\
\text{(resp. }\left(  g+f\right)  (a)  &  \leq\left(  g+f\right)
(b)\text{\ id est }f(a)-f(b)\leq g(b)-g(a)\text{).}%
\end{align*}
D'o\`{u} $\left\vert f(b)-f(a)\right\vert \leq g(b)-g(a)$.\medskip

\noindent A ce stade, nous avons obtenu le diagramme suivant\footnote{Dans ce
diagramme, ainsi que dans les diagrammes (\ref{FCD2008schemaEquivalences}) et
(\ref{FCD2008schemaEquivalences2}), les fl\`{e}ches unidirectionnelles sont
des implications math\'{e}matiques et les fl\`{e}ches bidirectionnelles des
\'{e}quivalences.}~:%
\begin{equation}%
\begin{tabular}
[c]{ccccccc}%
IAFG &  & $\longrightarrow$ &  & IAF &  & \\
& $\searrow\hspace{-0.15in}\nwarrow$ &  & $\nearrow\ast$ &  & $\searrow$ & \\
$\updownarrow$ &  & \textit{MAJA} &  & $\updownarrow$ &  & FCD\\
& $\nearrow\hspace{-0.15in}\swarrow$ &  & $\searrow\ast$ &  & $\nearrow\sharp$
& \\
SVD &  & $\longrightarrow$ &  & IAF' &  &
\end{tabular}
\label{FCD2008Diagramme1}%
\end{equation}
(Les implications $\left(  \ast\right)  $ sont imm\'{e}diates, en appliquant
la propri\'{e}t\'{e} MAJA \`{a} $f$ et $-f$. Pour l'implication $(\sharp)$, il
est en effet clair que si $f^{^{\prime}}$ est nulle, donc positive et
n\'{e}gative, $f$ est \`{a} la fois croissante et d\'{e}croissante, donc
constante.)\medskip

\noindent De mani\`{e}re claire, la propri\'{e}t\'{e} de majoration des
accroissements MAJA est centrale, puisqu'elle entra\^{\i}ne toutes les
autres.\ Nous allons nous int\'{e}resser aux d\'{e}monstrations possibles des
propri\'{e}t\'{e}s FCD, IAF et MAJA, en commen\c{c}ant ci-dessous par celle
d\'{e}coulant du th\'{e}or\`{e}me des accroissements finis.

\subsection{Les bases de la d\'{e}monstration de l'\'{e}galit\'{e} des
accroissements finis}

Rappelons l'\'{e}galit\'{e} des accroissements finis sous sa forme la plus \'{e}l\'{e}mentaire.

\begin{thm}
\label{FCD2008EgaliteACF}Soit $f:[a,b]\rightarrow\mathbb{R}$ continue sur
$[a,b]$ et d\'{e}rivable sur $\left]  a,b\right[  $. Alors il existe
$c\in\left]  a,b\right[  $ tel que $f(b)-f(a)=f^{\prime}(c)\left(  b-a\right)
$.
\end{thm}

\noindent Le th\'{e}or\`{e}me \ref{FCD2008EgaliteACF} entra\^{\i}ne
imm\'{e}diatement la caract\'{e}risation FCD, la caract\'{e}risation SVD, la
propri\'{e}t\'{e} MAJA et l'in\'{e}galit\'{e} classique des accroissements
finis.\medskip

\noindent La d\'{e}monstration la plus courante de l'\'{e}galit\'{e} des
accroissements finis consiste \`{a} la faire d\'{e}couler du th\'{e}or\`{e}me
de Rolle que nous rappelons ci-dessous.

\begin{thm}
\label{FCD2008theoremeRolle}Soit $g:[a,b]\rightarrow\mathbb{R}$ continue sur
$[a,b]$ et d\'{e}rivable sur $\left]  a,b\right[  $ telle que $g(a)=g(b)$%
.\ Alors il existe $c\in\left]  a,b\right[  $ tel que $g^{\prime}(c)=0$.
\end{thm}

\noindent Il suffit d'appliquer le th\'{e}or\`{e}me de Rolle \`{a} $g:x\mapsto
f(x)-\frac{f(b)-f(a)}{b-a}(x-a)$ pour obtenir le th\'{e}or\`{e}me
\ref{FCD2008EgaliteACF}. Notons, r\'{e}ciproquement, que l'application de
l'\'{e}galit\'{e} des accroissements finis avec $f:=g$ permet de v\'{e}rifier
que l'\'{e}galit\'{e} des accroissements finis et le th\'{e}or\`{e}me de Rolle
sont deux \'{e}nonc\'{e}s \'{e}quivalents.\ Focalisons nous un instant sur la
d\'{e}monstration de ce dernier. On utilise le plus souvent deux r\'{e}sultats
fondamentaux, l'un du cours sur les fonctions continues et l'autre du cours de
calcul diff\'{e}rentiel.\ On les rappelle ci dessous (th\'{e}or\`{e}mes
\ref{FCD2008theoremeFctContinue1} et \ref{FCD2008CNOrdre1}).

\begin{thm}
\label{FCD2008theoremeFctContinue1}(\textbf{Image d'un segment par une
application continue}) Soit $g:[a,b]\rightarrow\mathbb{R}$ continue.\ Alors
$g$ est born\'{e}e et atteint ses bornes.
\end{thm}

\begin{remarque}
\label{FCD2008ImageSegment}Nous ne retenons pas ici que $g\left(
[a,b]\right)  $ est un intervalle compact, cons\'{e}quence directe du
th\'{e}or\`{e}me des valeurs interm\'{e}diaires appliqu\'{e} \`{a} $g$, dans
la mesure o\`{u} cela n'est pas utile \`{a} la d\'{e}monstration du
th\'{e}or\`{e}me de Rolle.
\end{remarque}

\noindent De mani\`{e}re classique, dans un cours du d\'{e}but d'enseignement
sup\'{e}rieur, on d\'{e}montre le th\'{e}or\`{e}me
\ref{FCD2008theoremeFctContinue1} \`{a} l'aide de la propri\'{e}t\'{e} $(BW)$
de Bolzano-Weierstrass.\ On utilise $(BW)$ une premi\`{e}re fois pour
d\'{e}montrer que $g$ est born\'{e}e.\ Cette d\'{e}monstration se fait
facilement par l'absurde\footnote{On peut aussi utiliser la propri\'{e}t\'{e}
$(IE)$ des segments embo\^{\i}t\'{e}s~pour cette partie de la
d\'{e}monstration.}. On l'utilise une seconde fois pour d\'{e}montrer par un
raisonnement direct que $g$ atteint ses bornes.

\begin{thm}
\label{FCD2008CNOrdre1}(\textbf{Condition n\'{e}cessaire d'extremum du premier
ordre}) Soit $g:I\rightarrow\mathbb{R}$ d\'{e}rivable ($I$ intervalle
r\'{e}el).\ En tout extremum appartenant \`{a} l'int\'{e}rieur de $I$,
$g^{\prime}$ s'annule.
\end{thm}

\noindent Une fois le th\'{e}or\`{e}me \ref{FCD2008theoremeFctContinue1}
acquis, la d\'{e}monstration du th\'{e}or\`{e}me de Rolle est une simple
application du th\'{e}or\`{e}me \ref{FCD2008CNOrdre1}. Le cas o\`{u} $g$ est
constante \'{e}tant imm\'{e}diat, on se place dans la situation o\`{u} $g$ ne
l'est pas.\ La fonction $g$ poss\`{e}de alors un maximum et un minimum dont
l'un au moins est une valeur prise en $c\in\left]  a,b\right[  $. En
appliquant la condition n\'{e}cessaire du premier ordre pour l'existence d'un
extremun, on obtient $g^{\prime}(c)=0$.

\begin{remarque}
Notons que l'hypoth\`{e}se $g(a)=g(b)$ est exploit\'{e}e pour v\'{e}rifier que
la fonction $g$ poss\`{e}de un extremum pris \`{a} l'int\'{e}rieur de $\left[
a,b\right]  $.\ On retrouvera cette situation ult\'{e}rieurement.
\end{remarque}

\section{D\'{e}monstrations directes des caract\'{e}risations FCD et SVD par
un argument de connexit\'{e}\label{FCDPartieConnexite}}

Ces d\'{e}monstrations ont en commun \ un esprit qui est celui des classes
pr\'{e}paratoires des ann\'{e}es 1970.\ On en reprend l'id\'{e}e dans
l'ouvrage \textit{Math\'{e}matiques g\'{e}n\'{e}rales }(C. Pisot et
M.\ Zamanski, 1972). Leur principe servait aussi dans les ouvrages de
E.\ Ramis, dont la premi\`{e}re \'{e}dition serait de 1970\ ou 1971 (source:
http://publimath.irem.univ-mrs.fr/bibliocomp/\-M1U99048.htm, site consult\'{e}
le 25 janvier 2008) pour d\'{e}montrer l'in\'{e}galit\'{e} des accroissements
finis. Leur avantage est de s'\'{e}tendre facilement au cas d'une fonction $f$
\`{a} valeurs dans un espace norm\'{e}.

\subsection{Caract\'{e}risation FCD}

\begin{proposition}
\label{FCD2008DemoFCDConnexite}Soit $f:\left[  \alpha,\beta\right]
\rightarrow\mathbb{R}$ une fonction continue sur $\left[  \alpha,\beta\right]
$, admettant en tout point de $\left]  \alpha,\beta\right[  $ une
d\'{e}riv\'{e}e \`{a} droite nulle.\ Alors $f$ est constante.
\end{proposition}

\noindent\textbf{Preuve}.--~Nous restons proche du texte de C. Pisot et
M.\ Zamanski, cit\'{e} plus haut.\ Cependant, nous ne supposons pas la
d\'{e}rivabilit\'{e} \`{a} droite aux points $\alpha$ et $\beta$, ce qui
conduit \`{a} un petit artifice explicit\'{e} en fin de d\'{e}monstration.
Soit $a$, $b$ deux points de l'int\'{e}rieur de l'intervalle $I$ avec $a<b$.
Soit $\varepsilon$ strictement positif, donn\'{e} quelconque. Consid\'{e}rons
l'ensemble $X_{\varepsilon}$ des points $x$ de $[a,b]$ tels que%
\[
|f(x)-f(a)|/(x-a)\leq\varepsilon.
\]

\noindent Comme $f_{d}^{\prime}(a)$ existe et vaut z\'{e}ro, l'ensemble
$X_{\varepsilon}$ n'est pas vide et contient un intervalle d'extr\'{e}mit\'{e}
gauche $a$.\ Cet ensemble est major\'{e} par $b$.\ Comme selon la
propri\'{e}t\'{e} $\left(  BS\right)  $, toute partie de $\mathbb{R}$ non vide
major\'{e}e admet une borne sup\'{e}rieure, $X_{\varepsilon}$ poss\`{e}de une
borne sup\'{e}rieure $\xi$.\ Par une caract\'{e}risation de la borne
sup\'{e}rieure, il existe une suite de points $\left(  x_{n}\right)  _{n}$ de
$X_{\varepsilon}$ tendant vers $\xi$. Comme $|f(x_{n})-f(a)|\leq
\varepsilon\left(  x_{n}-a\right)  $ et $f$ est continue, on a $|f(\xi
)-f(a)|\leq\varepsilon\left(  \xi-a\right)  $.\ Donc $\xi\in X_{\varepsilon}$.

\noindent D\'{e}montrons maintenant par l'absurde que $\xi=b$. Supposons
$\xi<b$.\ Il existe alors $\varsigma\in\left]  \xi,b\right]  $ tel que
$|f(\varsigma)-f(\xi)|/\left(  \varsigma-\xi\right)  \leq\varepsilon$, en
utilisant $f_{d}^{\prime}(\xi)=0$ comme on a utilis\'{e} $f_{d}^{\prime}(a)=0$
ci-dessus. On a alors
\[
|f(\varsigma)-f(a)|\leq\underset{\leq\varepsilon\left(  \varsigma-\xi\right)
}{\underbrace{|f(\varsigma)-f(\xi)|}}+\underset{\leq\varepsilon\left(
\xi-a\right)  }{\underbrace{|f(\xi)-f(a)|}}\leq\varepsilon\left(
\varsigma-a\right)  .
\]
Ceci contredit le fait que $\xi$ soit la borne sup\'{e}rieure de l'ensemble
$X_{\varepsilon}$. On a donc $\xi=b$ et $X_{\varepsilon}=\left[  a,b\right]
$.\ Ainsi%
\[
|f(b)-f(a)|\leq\varepsilon(b-a).
\]
Comme on a l'in\'{e}galit\'{e} pr\'{e}c\'{e}dente pour $\varepsilon>0$
quelconque, on en d\'{e}duit $f(b)=f(a)$.\ En utilisant la continuit\'{e} de
$f$ sur $\left[  \alpha,\beta\right]  $, on en d\'{e}duit facilement que $f$
est constante sur $\left[  \alpha,\beta\right]  $.$\square$\medskip

\noindent En fait, moyennant des adaptations mineures, ce raisonnement montre
la fermeture de l'ensemble $X_{\varepsilon}$ (lorsqu'on montre que $\xi\in
X_{\varepsilon}$) et son ouverture (lorsqu'on montre que $\xi=b$).\ Ainsi,
$X_{\varepsilon}$ est un sous ensemble non vide, ferm\'{e} et ouvert dans
l'intervalle $I$, qui est connexe.\ C'est donc que $X_{\varepsilon}=I$.

\subsection{Caract\'{e}risation MAJA et SVD}

Il est tr\`{e}s tentant de s'inspirer de la d\'{e}monstration de la
proposition \ref{FCD2008DemoFCDConnexite} pour d\'{e}montrer la
caract\'{e}risation MAJA ou celle SVD sans r\'{e}f\'{e}rence du moins directe
\`{a} l'in\'{e}galit\'{e} des accroissements finis.

\begin{proposition}
\label{FCD2008DemoSVDConnexite}(Principe de Lagrange) Soit $f:\left[
\alpha,\beta\right]  \rightarrow\mathbb{R}$ une fonction continue sur $\left[
\alpha,\beta\right]  $, admettant en tout point de $\left]  \alpha
,\beta\right[  $ une d\'{e}riv\'{e}e \`{a} droite $f_{d}^{\prime}$ major\'{e}e
par un r\'{e}el $M\geq0$.\ Alors%
\[
f(b)-f(a)\leq M(b-a).
\]

\end{proposition}

\noindent Avant de donner les principales \'{e}tapes de la d\'{e}monstration,
notons que le cas $M=0$ donne la d\'{e}monstration de la propri\'{e}t\'{e}
SVD, pour $f_{d}^{\prime}$ n\'{e}gative.\medskip

\noindent Avec les m\^{e}mes notations que dans la d\'{e}monstration de la
proposition \ref{FCD2008DemoFCDConnexite}, on consid\`{e}re ici l'ensemble
$Y_{\varepsilon}$ des points $x$ de $[a,b]$ tels que%
\[
\left(  f(x)-f(a)\right)  /(x-a)\leq M+\varepsilon.
\]
L'ensemble $Y_{\varepsilon}$ est non vide (puisque $f_{d}^{\prime}(a)\leq M$)
major\'{e} par $b$ et sa borne sup\'{e}rieure $\xi$ lui appartient par
l'argument de continuit\'{e} d\'{e}j\`{a} employ\'{e}.\ De m\^{e}me, si on
suppose $\xi<b$, on trouve $\varsigma\in\left]  \xi,b\right]  $ tel que
$\left(  f(\varsigma)-f(\xi)\right)  /\left(  \varsigma-\xi\right)  \leq
M+\varepsilon$, en utilisant $f_{d}^{\prime}(\xi)\leq M$, comme on a
utilis\'{e} $f_{d}^{\prime}(\xi)=0$ ci-dessus. On a alors
\[
f(\varsigma)-f(a)=\underset{\leq\left(  M+\varepsilon\right)  \left(
\varsigma-\xi\right)  }{\underbrace{f(\varsigma)-f(\xi)}}+\underset
{\leq\left(  M+\varepsilon\right)  \left(  \xi-a\right)  }{\underbrace
{f(\xi)-f(a)}}\leq\left(  M+\varepsilon\right)  \left(  \varsigma-a\right)  .
\]
On contredit le fait que $\xi$ soit la borne sup\'{e}rieure. Il vient donc
$b=\xi$ et $f(b)-f(a)\leq\left(  M+\varepsilon\right)  \left(  b-a\right)  $
pour tout $\varepsilon>0$. La conclusion en d\'{e}coule.$\square$

\section{D\'{e}monstration des caract\'{e}risations FCD et SVD par un
processus de dichotomie\label{FCD2008DemoDicho}}

\subsection{Un lemme pr\'{e}paratoire}

Pour $f:I\rightarrow\mathbb{R}$ ($I$ intervalle r\'{e}el) et $\left(
x,y\right)  \in I^{2}$ $x\neq y$, posons
\[
P(x,y)=\frac{f\left(  y\right)  -f\left(  x\right)  }{y-x}.
\]
Dans un cadre g\'{e}om\'{e}trique, cette quantit\'{e}, la fonction "pente",
est interpr\'{e}table comme la pente de la corde de l'arc du graphe de $f$
d'extremit\'{e}s $\left(  x,f(x)\right)  $ et $\left(  y,f(y)\right)  $. On
remarque (simple calcul alg\'{e}brique) que $P$ est sym\'{e}trique en $x$ et
$y$ et qu'elle v\'{e}rifie, pour $\left(  x,y,a\right)  \in I^{3}$ deux \`{a}
deux disctincts,
\[
P(x,y)=\frac{y-a}{y-x}\frac{f(y)-f(a)}{y-a}+\frac{a-x}{y-x}\,\frac{f\left(
a\right)  -f\left(  x\right)  }{a-x}=\frac{a-x}{y-x}P(a,x)+\frac{y-a}%
{y-x}P(a,y).
\]

\begin{lemme}
\label{FCD2008LemmeAdad1}Soit $I$ un intervalle ouvert de $\mathbb{R}$ et
$f:I\rightarrow\mathbb{R}$ d\'{e}rivable en $a\in$ $I$. Alors%
\[
\lim_{x\rightarrow a^{-},y\rightarrow a^{+}}P(x,y)=f^{\prime}\left(  a\right)
.
\]
o\`{u} $\lim_{x\rightarrow a^{-},y\rightarrow a^{+}}$ veut dire limite pour
$x$ tendant vers $a$ par valeurs strictement inf\'{e}rieures et $y$ tendant
vers $a$ par valeurs strictement sup\'{e}rieures.
\end{lemme}

\noindent\textbf{Preuve}.--~Il suffit d'\'{e}crire
\begin{equation}
P(x,y))=\frac{a-x}{y-x}P(a,x)+\frac{y-a}{y-x}P(a,y)
\label{FCD2008RelationPentes}%
\end{equation}
avec%
\[
\underset{>0}{\underbrace{\frac{a-x}{y-x}}}+\underset{>0}{\underbrace
{\frac{y-a}{y-x}}}=1.
\]
\ La relation pr\'{e}c\'{e}dente est alors une relation barycentrique qui
montre que $P(y,x)$ appartient au segments d'extr\'{e}mit\'{e}s $P(x,a)$ et
$P(y,a)$.\ Par le th\'{e}or\`{e}me des gendarmes, on a $\lim_{x\rightarrow
a^{-},y\rightarrow a^{+}}P(y,x)=f(^{\prime}a)$.$\square$

\begin{remarque}
Le lecteur constatera que cette d\'{e}monstration utilise en fait une
propri\'{e}t\'{e} g\'{e}n\'{e}rale d'un triangle. On renvoie au
\ref{FCD2008SousSecPente} pour plus de commentaires \`{a} ce sujet.
\end{remarque}

\subsection{La d\'{e}monstration}

\noindent Pour une fonction $f$ continue sur un intervalle $I$ et
d\'{e}rivable sur l'int\'{e}rieur de $I$, la contrapos\'{e}e de la
propri\'{e}t\'{e}
\[
f^{\prime}\ \text{nulle sur l'int\'{e}rieur de }I\text{\ }\Rightarrow
\ f\text{\ constante sur }I
\]
est%
\[
f\text{\ non constante sur }I\Rightarrow f^{\prime}\ \text{non nulle sur
l'int\'{e}rieur de }I.
\]
On d\'{e}montrera donc le lemme suivant.

\begin{lemme}
\label{FCD2008LemmeFond1}(\textbf{Contrapos\'{e}e de la caract\'{e}risation
FCD}) Soit $f:\left[  \alpha,\beta\right]  \rightarrow\mathbb{R}$ continue sur
$\left[  \alpha,\beta\right]  $ et d\'{e}rivable sur $\left]  \alpha
,\beta\right[  $. On suppose qu'il existe $\left(  a,b\right)  \in\left[
\alpha,\beta\right]  ^{2}$ avec $a<b$ tels que $f(a)\neq f(b)$. Alors il
existe $c\in\left[  a,b\right]  $ tels que $f^{\prime}(c)\neq0$.
\end{lemme}

\begin{remarque}
Notons que l'on peut toujours supposer que $\left(  a,b\right)  \in\left]
\alpha,\beta\right[  ^{2}$, ce que nous ferons dans la suite. En effet, si par
exemple $a=\alpha$, le th\'{e}or\`{e}me des valeurs interm\'{e}diaires,
appliqu\'{e} \`{a} $f$ permet de trouver $a^{\prime}\in\left]  \alpha
,\beta\right]  $ tel que $f(a^{\prime})\neq f\left(  b\right)  $.
\end{remarque}

\noindent\textbf{Preuve}.--~On pose
\[
a_{0}=a\ ,\ \ b_{0}=b\ ;\ \ d=\left\vert f(b)-f(a)\right\vert .
\]
\textsf{Etape 1}.--~On construit par r\'{e}currence deux suites adjacentes
$\left(  a_{n}\right)  _{n}$ et $\left(  b_{n}\right)  _{n}$ telles que, pour
tout $n\in\mathbb{N}$,%
\begin{equation}
\left\vert f\left(  b_{n}\right)  -f\left(  a_{n}\right)  \right\vert
\geq\frac{d}{2^{n}}\ ,\ \ b_{n}-a_{n}=\frac{b-a}{2^{n}}.
\label{CroissanceSansTAF1}%
\end{equation}
\textit{Initialisation}.--~$a_{0}$ et $b_{0}$ satisfont les relations
(\ref{CroissanceSansTAF1}) pour $n=0$.

\noindent\textit{H\'{e}r\'{e}dit\'{e}}.--~Soit $n\geq0$.\ On suppose
construits $a_{n}$ et $b_{n}$ satisfaisant les relations
(\ref{CroissanceSansTAF1}). Posons $c_{n}=\left(  a_{n}+b_{n}\right)
/2$.\ Comme $\left\vert f\left(  b_{n}\right)  -f\left(  a_{n}\right)
\right\vert \geq d/2^{n}$, on a n\'{e}cessairement $\left\vert f\left(
c_{n}\right)  -f\left(  a_{n}\right)  \right\vert \geq d/2^{n+1}$ ou
$\left\vert f\left(  b_{n}\right)  -f\left(  c_{n}\right)  \right\vert \geq
d/2^{n+1}$.\ Si $\left\vert f\left(  c_{n}\right)  -f\left(  a_{n}\right)
\right\vert \geq d/2^{n+1}$, on pose $a_{n+1}=a_{n}$ et $b_{n+1}=c_{n}%
$.\ Sinon, on pose $a_{n+1}=c_{n}$ et $b_{n+1}=b_{n}$. De mani\`{e}re claire,
$a_{n+1}$ et $b_{n+1}$ satisfont les relations (\ref{CroissanceSansTAF1}) et
la suite $\left(  a_{n}\right)  _{n}$ est croissante, la suite $\left(
b_{n}\right)  _{n}$ d\'{e}croissante.\medskip

\noindent Les suite $\left(  a_{n}\right)  _{n}$ et $\left(  b_{n}\right)
_{n}$ convergent, comme suites adjacentes, vers la m\^{e}me limite
$c\in\left[  a,b\right]  $. C'est, en fait, la propri\'{e}t\'{e} $(SE)$\emph{
}des segments embo\^{\i}t\'{e}s qui est utilis\'{e}e.\emph{ }\medskip

\noindent\textsf{Etape 2}.--~Supposons l'une des suites $\left(  a_{n}\right)
_{n}$ ou $\left(  b_{n}\right)  _{n}$ constante \`{a} partir d'un certain rang
(ou stationnaire). Par exemple, supposons qu'il s'agisse de la suite $\left(
a_{n}\right)  _{n}$.\ Il existe $m\in\mathbb{N}$ tel que $\forall n\geq
m,\ a_{n}=a_{m}=c$. On a alors, pout tout $n\geq m$, $a_{n}=c<b_{n}$.\ Alors,
on peut consid\'{e}rer le quotient
\[
P\left(  b_{n},c\right)  =\frac{f\left(  b_{n}\right)  -f\left(  c\right)
}{b_{n}-c}%
\]
avec lim$_{n\rightarrow+\infty}P_{n}\left(  b_{n},c\right)  =f^{\prime}%
(c)$.\ Compte tenu des relations (\ref{CroissanceSansTAF1}), on a
\[
\left\vert P_{n}\left(  b_{n},c\right)  \right\vert \geq\frac{d}{2^{n}}%
\frac{2^{n}}{b-a}=\frac{d}{b-a}>0.
\]
Par prolongement d'\'{e}galit\'{e}, il vient$\left\vert f^{\prime
}(c)\right\vert \geq d/\left(  b-a\right)  >0$.\medskip

\noindent\textsf{Etape 3}.--~On se place dans le cas ou les suites $\left(
a_{n}\right)  _{n}$ et $\left(  b_{n}\right)  _{n}$ ne sont pas stationnaires.
Rappelons que l'on a, pour tout $n\in\mathbb{N}$, $a_{n}\leq c$.\ Avec
l'hypoth\`{e}se de non stationarit\'{e} de la suite $\left(  a_{n}\right)
_{n}$, on a n\'{e}cessairement~: $\forall n\in\mathbb{N},\ a_{n}<c$. (En
effet, si il existait $m\in\mathbb{N}$ tel que $a_{m}=c$, on aurait, par
croissance de la suite $\left(  a_{n}\right)  _{n}$, la propri\'{e}t\'{e}~:
$\forall n\geq m,\ a_{n}=c$, ce qui nous ram\'{e}nerait \`{a} l'\'{e}tape 2.)
De m\^{e}me, on a~: $\forall n\in\mathbb{N},\ c<b_{n}$.\ On peut alors
appliquer la propri\'{e}t\'{e} d\'{e}montr\'{e}e dans le lemme%
\[
P\left(  b_{n},a_{n}\right)  =\frac{f\left(  b_{n}\right)  -f\left(
a_{n}\right)  }{b_{n}-a_{n}}\overset{n\rightarrow+\infty}{\longrightarrow
}f^{\prime}(c),
\]
pour conclure comme dans l'\'{e}tape 2, puisque
\[
\left\vert P\left(  b_{n},a_{n}\right)  \right\vert =|\frac{f\left(
b_{n}\right)  -f\left(  a_{n}\right)  }{b_{n}-a_{n}}|\geq\frac{d}%
{b-a}>0.\square
\]

\begin{remarque}
Le cas d'une suite $\left(  a_{n}\right)  _{n}$ ou $\left(  b_{n}\right)
_{n}$ stationnaire peut survenir.\ Par exemple, en prenant $a=0$, $b=1$ et
$f(x)=x$, on voit qu'un choix possible est de prendre la suite $\left(
a_{n}\right)  _{n}$ nulle et la suite $\left(  b_{n}\right)  _{n}$ \'{e}gale
\`{a} $\left(  1/2^{n}\right)  _{n\geq0}$.
\end{remarque}

\begin{question}
Peut-on toujours construire $\left(  a_{n}\right)  _{n}$ et $\left(
b_{n}\right)  _{n}$ non stationnaires~?
\end{question}

\noindent Notons que le lemme \ref{FCD2008LemmeFond1} poss\`{e}de son exacte
contrepartie pour le principe de Lagrange.

\begin{lemme}
\label{FCD2008LemmeFondSVD}(\textbf{Contrapos\'{e}e du principe de Lagrange})
Soit $f:\left[  \alpha,\beta\right]  \rightarrow\mathbb{R}$ continue sur
$\left[  \alpha,\beta\right]  $ et d\'{e}rivable sur $\left]  \alpha
,\beta\right[  $. On suppose qu'il existe $\left(  a,b\right)  \in\left[
\alpha,\beta\right]  ^{2}$ avec $a<b$ tels que $f(a)<f(b)$ (resp.
$f(a)>f(b)$). Alors il existe $c\in\left[  a,b\right]  $ tels que $f^{\prime
}(c)>0$ (resp. $f^{\prime}(c)<0$).
\end{lemme}

\noindent Le lecteur s'assurera de l'\'{e}quivalence du lemme
\ref{FCD2008LemmeFondSVD} et du principe de Lagrange soit par un raisonnement
par l'absurde soit en \'{e}crivant le lemme \ref{FCD2008LemmeFondSVD} et le
principe de Lagrange sous forme de propositions logiques formalis\'{e}es. La
d\'{e}monstration du lemme \ref{FCD2008LemmeFondSVD} par dichotomie reprend
exactement celle du lemme, en la simplifiant m\^{e}me l\'{e}g\'{e}rement, dans
la mesure o\`{u} l'on supprime les valeurs absolues. (On pose $d=f(b)-f(a).$)
\medskip

\noindent Plus int\'{e}ressant, cette d\'{e}monstration fournit \'{e}galement
une preuve de la propri\'{e}t\'{e} de majoration des accroissements ou de
l'in\'{e}galit\'{e} des accroissements finis classique par dichotomie. Par
exemple, sous les hypoth\`{e}ses de la Proposition
\ref{FCD2008LemmeACFClassique}, supposons $\left\vert f(b)-f(a)\right\vert
>k(b-a)$.\ On pose alors $d=f(b)-f(a)=k^{\prime}(b-a)$ avec $k^{\prime}>k$
dans la preuve du lemme \ref{FCD2008LemmeFond1}.\ On obtient alors l'existence
de $c\in\left[  a,b\right]  $ tel que $\left\vert f^{\prime}(c)\right\vert
\geq k^{\prime}>k$, contredisant l'hypoth\`{e}se.\ 

\section{Compl\'{e}ments dans le cadre de l'analyse r\'{e}elle
\'{e}l\'{e}mentaire\label{FCD2008ComMath}}

La d\'{e}monstration par dichotomie et par le th\'{e}or\`{e}me des
accroissements finis montre le r\^{o}le central de la fonction "pente",
fonction sur laquelle non donnons quelques compl\'{e}ments ci-dessous.\ Ces
d\'{e}monstrations permettent \'{e}galement de se r\'{e}interroger sur un
c\'{e}l\`{e}bre th\'{e}or\`{e}me de Darboux, affirmant que les
d\'{e}riv\'{e}es poss\`{e}dent la propri\'{e}t\'{e} des valeurs
interm\'{e}diaires. Enfin, nous pr\'{e}sentons quelques ouvertures
math\'{e}matiques autour de la caract\'{e}risation FCD.

\subsection{Pente et d\'{e}riv\'{e}e\label{FCD2008SousSecPente}}

Le lemme \ref{FCD2008LemmeAdad1} est compl\'{e}t\'{e} par la proposition suivante.

\begin{proposition}
\label{FCD2008PropPenteDerivee}Soit $I$ un intervalle ouvert de $\mathbb{R}$,
$f:I\rightarrow\mathbb{R}$ et $a\in$ $I$. La fonction $f$ est d\'{e}rivable en
$a$ de nombre d\'{e}riv\'{e}e $l$ si, et seulement si,%
\begin{equation}
\lim_{x\rightarrow a^{-},y\rightarrow a^{+}}P(x,y)=l. \label{FCD2008CNDPente}%
\end{equation}

\end{proposition}

\noindent On note que si l'on remplace, dans la proposition
\ref{FCD2008PropPenteDerivee}, la notion de limite consid\'{e}r\'{e}e par la
notion de limite pour $x$ (resp. $y$) tendant vers $a$ par valeurs
inf\'{e}rieures (resp. sup\'{e}rieures) ou \'{e}gales \`{a} $a$ et $x\neq y$,
la d\'{e}rivabilit\'{e} en $a$ devient \'{e}vidente puisque $f$ sera alors
d\'{e}rivable \`{a} droite et \`{a} gauche en $a$, avec \'{e}galit\'{e} de ces
deux d\'{e}riv\'{e}es unilat\'{e}rales \cite{DoubkhanSifre}. En revenant \`{a}
la notion de limite intialement consid\'{e}r\'{e}e, montrons que l'existence
de $l$ rendant vraie la relation (\ref{FCD2008CNDPente}) entra\^{\i}ne que
$f^{\prime}(a)=l$. La relation (\ref{FCD2008CNDPente}) peut se
r\'{e}\'{e}crire sous la forme%
\[
\lim_{h\rightarrow0^{+},k\rightarrow0^{+}}P(a-h,a+k)=l
\]
et la relation (\ref{FCD2008RelationPentes}) sous la forme%
\[
P(a-h,a+k)=\frac{h}{h+k}P(a-h,a)+\frac{k}{h+k}P(a,a+k).
\]
Ainsi%
\[
\frac{h+k}{h}P(a-h,a+k)=P(a-h,a)+\frac{k}{h}P(a,a+k)=P(a-h,a)+\frac{1}%
{h}\left(  f(a+k)-f(a)\right)  .
\]
La continuit\'{e} de $f$ en $a$ permet de trouver, pour tout $h>0$, un
r\'{e}el $k(h)$ tel que $\lim_{h\rightarrow0^{+}}k(h)=0$ et $\left\vert
f(a+k(h))-f(a)\right\vert <h^{2}$. De plus, on peut supposer que
$\lim_{h\rightarrow0^{+}}k(h)/h=0$. Alors
\[
\lim_{h\rightarrow0^{+}}\left(  1+k(h)/h\right)  =1\ ,\ \lim_{h\rightarrow
0^{+}}P(a-h,a+k(h))=l\ ,\ \ \lim_{h\rightarrow0^{+}}\frac{1}{h}\left(
f(a+k(h))-f(a)\right)  =0.
\]
D'o\`{u} $\lim_{h\rightarrow0^{+}}P(a-h,a)=l$.\ Un raisonnement analogue
montre que $\lim_{k\rightarrow0^{+}}P(a,a+k)=l.$

\begin{remarques}
\label{FCD2008PenteDeriveeContinue}~\newline$\left(  i\right)  $~Dans la
proposition \ref{FCD2008PropPenteDerivee}, la continuit\'{e} au point $a$ doit
\^{e}tre suppos\'{e}e, pour \'{e}viter le cas d'une fonction mal d\'{e}finie
au point $a$. On prend $I=\mathbb{R}$, $f$ d\'{e}finie par $f(0)=0$ et
$f(x)=c$, constante non nulle, pour tout $x\in\mathbb{R}^{\ast}$.\ On a
$\lim_{h\rightarrow0^{+},k\rightarrow0^{+}}P(a-h,a+k)=0$, tandis que $f$ n'est
pas d\'{e}rivable en $0$ puisque non continue en ce point.\newline$\left(
ii\right)  $~La d\'{e}monstration ci-dessus poss\`{e}de une interpr\'{e}tation
heuristique (formalisable, \'{e}ventuellement, dans le langage de la
Th\'{e}orie Relative des Ensembles Internes de Y.\ P\'{e}raire \cite{Peraire1}%
) assez imm\'{e}diate~: le r\'{e}el $h$ \'{e}tant un infiniment petit du
premier ordre, $k$ un infiniment petit du second ordre, la pente $P(a-h,a+k)$
est infiniment proche au premier ordre de $l$.\ De plus, la diff\'{e}rence
entre $P(a-h,a+k)$ et $P(a-h,a)$ est un infiniment petit du premier ordre, car
la diff\'{e}rence entre $f(a+k)$ et $f(a)$ est un infiniment petit du
deuxi\`{e}me ordre et $h$ un infiniment petit du premier ordre. Ainsi
$P(a-h,a)$ est infiniment proche au premier ordre de $l$.\ (Le lecteur pourra
s'appuyer sur un dessin.)
\end{remarques}

\noindent La proposition \ref{FCD2008PropPenteDerivee} peut introduire une
mani\`{e}re de se requestionner sur le rapport entre d\'{e}rivabilit\'{e} et
stricte d\'{e}rivabilit\'{e} (voir, par exemple, \cite{MichelP1}).

\begin{defin}
\label{FCD2008DefStricteDeriv}Soit $I$ un intervalle r\'{e}el ouvert non
trivial, $f:I\rightarrow\mathbb{R}$, $a\in I$.\ On dit que $f$ est strictement
d\'{e}rivable en $a$ si $\lim_{(x,y)\rightarrow(a,a),x\neq y}P(x,y)$ existe.\ 
\end{defin}

\noindent Ce nombre, lorsqu'il existe est appel\'{e} la stricte
d\'{e}riv\'{e}e de $f$ au point $a$.\ On dit que $f$ est \emph{strictement
d\'{e}rivable sur }$I$ si $f$ est strictement d\'{e}rivable en tout point de
$I$. De mani\`{e}re claire, si $f$ est strictement d\'{e}rivable en $a\in I$,
$f$ est d\'{e}rivable en $a$. En revanche, $f$ peut \^{e}tre d\'{e}rivable en
$a$ sans \^{e}tre strictement d\'{e}rivable en $a$.

\begin{exem}
\label{FCD2008ExempleDerNSD1}On utilise la c\'{e}l\`{e}bre famille de
fonctions%
\[
f_{p,q}:\mathbb{R\rightarrow R}\ ,\ \ f(x)=x^{p}\sin(1/x^{q})\ \text{si }%
x\neq0,\ \ f(x)=0\ \text{si }x=0.
\]
avec, ici $p=2$,$\ q=1$. En posant $f_{2,1}=f$, on a
\[
f(x)=f\left(  x\right)  -f(0)=x\mathrm{o}\left(  1\right)  \ \mathrm{pour}%
\ x\rightarrow0.
\]
Ainsi, la fonction $f$ est d\'{e}rivable en $0$ avec $f^{\prime}\left(
0\right)  =0$. Cependant consid\'{e}rons les suites%
\[
\forall n\in\mathbb{N},\,x_{n}=\left(  \pi/2+\left(  2n+1\right)  \pi\right)
^{-1};\;y_{n}=\left(  \pi/2+2n\pi\right)  ^{-1}.
\]
On a%
\[
\frac{f\left(  y_{n}\right)  -f\left(  x_{n}\right)  }{y_{n}-x_{n}}%
\overset{n\rightarrow+\infty}{\longrightarrow}\frac{2}{\pi}\neq0=f^{\prime
}\left(  0\right)  .
\]
Ainsi $f$ n'est pas stritement d\'{e}rivable en $0$.
\end{exem}

\begin{remarque}
\label{FCD2008RemarqueLemmeDer}Cet exemple montre \'{e}galement qu'il est
essentiel dans le lemme \ref{FCD2008LemmeAdad1} que la limite soit prise pour
$x$ tendant vers $a$ par valeurs (strictement) inf\'{e}rieures et pour $y$
tendant vers $a$ par valeurs (strictement) sup\'{e}rieures.
\end{remarque}

\noindent De fait, les fonctions strictement d\'{e}rivables de la variable et
\`{a} valeurs r\'{e}elles sont bien connues d'apr\`{e}s le\emph{
th\'{e}or\`{e}me de la pente} rappel\'{e} ci-dessous.

\begin{thm}
\label{FCD2008PropPropSD1}(\textbf{Th\'{e}or\`{e}me de la pente}) Soit $I$ un
intervalle r\'{e}el ouvert non trivial et $f:I\rightarrow\mathbb{R}$. Sont
\'{e}quivalentes~:\newline$\left(  1\right)  $~Si $f$ est strictement
d\'{e}rivable sur $I$.\newline$\left(  2\right)  $~$f$ est de classe
$\mathrm{C}^{1}$ sur $I$.
\end{thm}

\noindent\textbf{Preuve}.--~Le $\left(  2\right)  \Rightarrow\left(  1\right)
$ peut se d\'{e}duire directement du th\'{e}or\`{e}me des accroissements
finis. Soit, en effet, $a\in I$. On \'{e}crit, pour tout $\left(  x,y\right)
\in I^{2}$ avec $x\neq y$,
\[
f(y)-f(x)=f^{\prime}(c_{x,y})(y-x)\ \ \text{avec }c_{x,y}\text{ compris entre
}x\text{ et }y.
\]
En utilisant la continuit\'{e} de $f^{\prime}$ en $a$, on a imm\'{e}diatement
$\lim_{\left(  x,y\right)  \rightarrow(a,a),x\neq y}P(x,y)=f^{\prime}%
(a)$.\ \smallskip

\noindent Pour le $\left(  1\right)  \Rightarrow\left(  2\right)  $, nous
allons proposer une d\'{e}monstration \'{e}l\'{e}mentaire bas\'{e}e sur les
d\'{e}finitions. Soit $a\in I$. Comme $f$ est suppos\'{e}e strictement
d\'{e}rivable sur $I$, pour tout $\varepsilon>0$, il existe $\eta_{a}>0$ tel
que $\left]  a-\eta_{a},a+\eta_{a}\right[  \subset I$ et
\begin{equation}
0<\left\vert x-a\right\vert <\eta_{a}\ \mathrm{et\ }0<\left\vert
y-a\right\vert <\eta_{a}\Rightarrow\left\vert P(x,y)-f^{\prime}(a)\right\vert
<\varepsilon/2.\ \label{FCD2008StricteDerivThP}%
\end{equation}
Soit maintenant $b\in I$ tel que $\left\vert b-a\right\vert <\eta=\eta_{a}%
/2$.\ Soit $\eta_{b}>0$ un r\'{e}el rendant vrai (\ref{FCD2008StricteDerivThP}%
) pour $b$ au lieu de $a$.\ Consid\'{e}rons $\left(  x,y\right)  \in I^{2}$
tels que $x\neq y$ et $\left\vert b-x\right\vert <\min\left(  \eta_{a}%
/2,\eta_{b}\right)  $, $\left\vert b-y\right\vert <\min\left(  \eta_{a}%
/2,\eta_{b}\right)  $.\ (Le lecteur s'assurera de l'existence d'un tel couple
$\left(  x,y\right)  $.) On a $\left\vert P(x,y)-f^{\prime}(b)\right\vert
<\varepsilon/2$ par d\'{e}finition de $\eta_{b}$.\ Comme $\left\vert
x-a\right\vert <\left\vert x-b\right\vert +\left\vert b-a\right\vert <\eta
_{a}$ et de m\^{e}me $\left\vert y-a\right\vert <\eta_{a}$, on a \'{e}galement
$\left\vert P(x,y)-f^{\prime}(a)\right\vert <\varepsilon/2$.\ D'o\`{u}
$\left\vert f^{\prime}(b)-f^{\prime}(a)\right\vert <\varepsilon$. Ainsi,
\[
\forall\varepsilon>0,\ \exists\eta>0,\ \forall b\in I,\ \ \left\vert
b-a\right\vert <\eta\Rightarrow\left\vert f^{\prime}(b)-f^{\prime
}(a)\right\vert <\varepsilon.\square
\]

\begin{remarque}
$f$ peut \^{e}tre strictement d\'{e}rivable en $a$ sans \^{e}tre d\'{e}rivable
sur un intervalle voisinage de $a$.
\end{remarque}

\begin{proposition}
Si $f$ est strictement d\'{e}rivable en $a\in I$, alors $f$ est lipschitzienne
au voisinage de $a$.
\end{proposition}

\noindent Nous renvoyons le lecteur int\'{e}ress\'{e} aux ouvrages classiques,
par exemple \cite{FClarke1}, pour plus de d\'{e}tails sur la stricte
d\'{e}rivabilit\'{e}. Pour le lien entre la notion de d\'{e}riv\'{e}e et la
tangente au graphe de la fonction (que nous avons utilis\'{e} implicitement
pour introduire la fonction des pentes et pour l'argument g\'{e}om\'{e}trique
\`{a} l'appui de la d\'{e}monstration du lemme \ref{FCD2008LemmeAdad1}), nous
renvoyons \`{a} la th\`{e}se de E.\ Rouy \cite{RouyTh}.

\subsection{Le lemme \ref{FCD2008LemmeFond1} et le th\'{e}or\`{e}me de
Darboux}

Rappelons le r\'{e}sultat suivant qui donne le second exemple de fonctions,
apr\`{e}s les fonctions continues, poss\'{e}dant la propri\'{e}t\'{e} de la
valeur interm\'{e}diaire.

\begin{thm}
\label{ddthmDbx}Soit $f:\left[  \alpha,\beta\right]  \rightarrow\mathbb{R}$
une fonction continue sur $\left[  \alpha,\beta\right]  $ et d\'{e}rivable sur
$\left]  \alpha,\beta\right[  $ L'ensemble $f^{\prime}(\left]  \alpha
,\beta\right[  )$ est un intervalle.
\end{thm}

\noindent Ce th\'{e}or\`{e}me est souvent pr\'{e}sent\'{e} comme une
application de l'\'{e}galit\'{e} des accroissements finis.\ Pour s'en
convaincre, donnons les id\'{e}es principales d'une d\'{e}monstration tr\`{e}s
\'{e}l\'{e}mentaire et se limitant \`{a} consid\'{e}rer des fonctions de la
variable r\'{e}elle. Il s'agit de montrer que, pour tout $(a,b)\in\left]
\alpha,\beta\right[  ^{2}$ avec $a<b$, le segment d'extr\'{e}mit\'{e}s
$f^{\prime}(a)$ et $f^{\prime}(b)$ est inclus dans $f^{\prime}(\left]
\alpha,\beta\right[  )$.\ On introduit alors la fonction $\varphi
:[a,b]\rightarrow\mathbb{R}$ (\emph{resp}. $\psi:[a,b]\rightarrow\mathbb{R}$)
d\'{e}finie%
\[
\varphi(a)=f^{\prime}(a)\;;\;\forall x\in\left]  a,b\right]  ,\;\varphi
(x)=P(a,x)\ \text{(\emph{resp}. }\psi(b)=f^{\prime}(b)\;;\;\forall x\in\left[
a,b\right[  ,\;\psi(x)=P(x,b)\text{).}%
\]
Le th\'{e}or\`{e}me des accroissements finis permet d'\'{e}crire%
\begin{equation}
\forall x\in\left]  a,b\right]  ,\ \exists c_{x}\in\left]  a,x\right[
:\ \ \varphi(x)=f^{\prime}(c_{x}). \label{FCD2008DemoDarboux1}%
\end{equation}
La propri\'{e}t\'{e} (\ref{FCD2008DemoDarboux1}) entra\^{\i}ne que
$\varphi([a,b])\subset f^{\prime}(\left]  \alpha,\beta\right[  )$. De plus, la
d\'{e}rivabilit\'{e} en $a$ de $f$ entra\^{\i}ne que $\lim_{x\rightarrow
a^{+}}\varphi(x)=f^{\prime}(a)=\varphi(a)$ ce qui \'{e}tait la seule
difficult\'{e} pour v\'{e}rifier que $\varphi$ est continue sur $[a,b]$.\ De
m\^{e}me on montre la continuit\'{e} de $\psi$ sur $[a,b]$ et l'inclusion
$\psi([a,b])\subset f^{\prime}(\left]  \alpha,\beta\right[  )$. Alors,
$\varphi([a,b])$ et $\psi([a,b])$ sont des intervalles. Comme $\varphi
(b)=\psi(a)$, on a $\varphi([a,b])\cap\psi([a,b])\neq\varnothing$ et la
r\'{e}union $\varphi([a,b])\cup\psi([a,b])$ est un intervalle inclus dans
$f^{\prime}(\left]  \alpha,\beta\right[  )$. De plus, elle contient le segment
d'extr\'{e}mit\'{e}s $f^{\prime}(a)$ et $f^{\prime}(b)$, ce qui entra\^{\i}ne
la conclusion\footnote{Une d\'{e}monstration, d'ailleurs plus courte, fait
d\'{e}couler cette propri\'{e}t\'{e} de l'ensemble form\'{e} des
th\'{e}or\`{e}mes \ref{FCD2008theoremeFctContinue1} et \ref{FCD2008CNOrdre1}%
.\ Cette d\'{e}monstration cour-circuite, en quelque sorte, le passage par le
th\'{e}or\`{e}me de Rolle et l'\'{e}galit\'{e} des accroissements finis. Voir
\cite{DoubkhanSifre}, qui renvoie \`{a} W. Rudin \cite{Rudin2} pour cette
d\'{e}monstration.}.\medskip

\noindent En fait, le th\'{e}or\`{e}me de Darboux, associ\'{e} au lemme
\ref{FCD2008LemmeFond1} (ou \`{a} la caract\'{e}risation SVD), est
\'{e}quivalent \`{a} l'\'{e}galit\'{e} des accroissements finis comme le
montre le raisonnement suivant. Soit $f:\left[  \alpha,\beta\right]
\rightarrow\mathbb{R}$ continue sur $\left[  \alpha,\beta\right]  $ et
d\'{e}rivable sur $\left]  \alpha,\beta\right[  $ telle que $f^{\prime}$
(d\'{e}finie sur $\left]  \alpha,\beta\right[  $) poss\`{e}de la
propri\'{e}t\'{e} de la valeur interm\'{e}diaire. Supposons de plus que
$f(\alpha)=f(\beta)$, hypoth\`{e}se du th\'{e}or\`{e}me de Rolle. Si $f$ est
constante sur $\left[  \alpha,\beta\right]  $, on sait que $f^{\prime}$ est
nulle et la conclusion du th\'{e}or\`{e}me de Rolle est satisfaite.\ Si $f$
n'est pas constante sur $\left[  \alpha,\beta\right]  $, il existe $\gamma
\in\left]  \alpha,\beta\right[  $ tel que $f(\alpha)\neq f(\gamma
)$.\ Supposons par exemple $f(\gamma)>f(\alpha)$. Il existe alors $a\in\left]
\alpha,\gamma\right[  $ tel que $f\left(  a\right)  <f\left(  \gamma\right)
$, d'apr\`{e}s le th\'{e}or\`{e}me des valeurs interm\'{e}diaires appliqu\'{e}
\`{a} $f$.\ De m\^{e}me, il existe $b\in\left]  \gamma,\beta\right[  $ tel que
$f\left(  \gamma\right)  >f(b)$. Le lemme \ref{FCD2008LemmeFond1} montre alors
qu'il existe $c\in\left[  a,\gamma\right]  $ tel que $f^{\prime}(c)>0$ et
$d\in\left[  \gamma,b\right]  $ tel que $f^{\prime}(d)<0$.\ D'apr\`{e}s la
propri\'{e}t\'{e} de valeur interm\'{e}diaire appliqu\'{e}e \`{a} $f^{\prime}%
$, il existe $d\in\left]  a,b\right[  $ tel que $f^{\prime}(d)=0$.\ On passe
alors du th\'{e}or\`{e}me de Rolle au th\'{e}or\`{e}me des accroissements
finis de mani\`{e}re classique.

\begin{remarque}
Le petit d\'{e}tour par les points $a$ et $b$ vient du fait que $f$ n'est pas
suppos\'{e}e d\'{e}rivable en $\alpha$ et $\beta$, alors que le lemme
\ref{FCD2008LemmeFond1} n\'{e}cessite a priori la d\'{e}rivabilit\'{e} aux
extr\'{e}mit\'{e}s du segment d'application.
\end{remarque}

\noindent Ainsi, on a les \'{e}quivalences~:%
\begin{equation}%
\begin{tabular}
[c]{ccc}%
Th\'{e}or\`{e}me de Rolle & $\longleftrightarrow$ & Th\'{e}or\`{e}me des
accroissements finis\\
$\hspace{0.9cm}\searrow\hspace{-0.15in}\nwarrow$ &  & $\nearrow\hspace
{-0.15in}\swarrow\hspace{1.1cm}$\\
\multicolumn{3}{c}{Th\'{e}or\`{e}me de Darboux\ \textit{et}%
\ Caract\'{e}risation SVD.$\hspace{1.15cm}$}%
\end{tabular}
\label{FCD2008schemaEquivalences}%
\end{equation}
Les caract\'{e}risations FCD et SVD s'ins\`{e}rent de la mani\`{e}re suivante
dans ce diagramme\footnote{L'implication $(\ast)$ d\'{e}coule d'une
application directe de l'\'{e}galit\'{e} des accroisssements finis au couple
$(a,b)$ tel que $f(a)<f(b)$ (resp. $f(a)>f(b)$.) Voir l'\'{e}nonc\'{e} du
lemme \ref{FCD2008LemmeFond1}.}~:%

\begin{equation}%
\begin{tabular}
[c]{lcrcl}
&  &
\begin{tabular}
[c]{c}%
In\'{e}galit\'{e}\\
accroissements finis
\end{tabular}
& $\longrightarrow$ & \textit{Caract\'{e}risation FCD}\\
\multicolumn{1}{r}{} &  & \multicolumn{1}{c}{$\uparrow$} & $_{\ast}%
\nwarrow\hspace{-0.16in}\nearrow$ & \multicolumn{1}{c}{$\uparrow$}\\
\multicolumn{1}{r}{Th\'{e}or\`{e}me de Rolle} & $\longleftrightarrow$ &
\begin{tabular}
[c]{c}%
Th\'{e}or\`{e}me des\\
accroissements finis
\end{tabular}
& $\longrightarrow$ & Caract\'{e}risation SVD\\
\multicolumn{1}{c}{$\updownarrow$} &  & \multicolumn{1}{c}{$\updownarrow$} &
$\nearrow$ & \\
\multicolumn{3}{r}{$\hspace{0.5cm}$Th\'{e}or\`{e}me de Darboux\ \textit{et}%
\ Caract\'{e}risation SVD.} &  &
\end{tabular}
\label{FCD2008schemaEquivalences2}%
\end{equation}
On a mis en valeur les techniques de d\'{e}monstration de la
caract\'{e}risation FCD, par les propri\'{e}t\'{e}s des accroissements finis
et par le principe de Lagrange (il manque ici les d\'{e}monstrations par le
raisonnement de connexit\'{e} et par le processus de Dichotomie.)

\begin{question}
La caract\'{e}risation FCD, associ\'{e}e au th\'{e}or\`{e}me de Darboux,
entra\^{\i}ne-telle une des propri\'{e}t\'{e}s \'{e}quivalentes du sch\'{e}ma
(\ref{FCD2008schemaEquivalences})~?
\end{question}

\subsection{La non trivialit\'{e} de la caract\'{e}risation FCD dans le cadre
des fonctions de la variable r\'{e}elle\label{FCD2008SPAnFonc}}

On rappelle ici que la conclusion des caract\'{e}risations FCD et SVD restent
valables sous des hypoth\`{e}ses l\'{e}g\'{e}rement plus faibles.\ Pour des
raisons d'homog\'{e}n\'{e}it\'{e} de l'expos\'{e}, on suppose donn\'{e}e
$f:\left[  a,b\right]  \rightarrow\mathbb{R}$ continue.

\begin{thm}
\label{FCD2008FCDSVDgeneralisee}Si le nombre d\'{e}riv\'{e}e $f^{\prime
}\left(  x\right)  $ existe sauf, \'{e}ventuellement, pour les points d'un
ensemble d\'{e}nombrable et s'il est nul (resp. positif) lorsqu'il existe,
alors $f$ est constante (resp.\ croissante.)
\end{thm}

\noindent Pour la d\'{e}monstration de ce r\'{e}sultat, nous renvoyons le
lecteur au cours d'analyse de J.\ Dieudonn\'{e} \cite{Dieudonne1}, tout en
mentionnant qu'elle constitue, en fait, un raffinement des d\'{e}monstrations
pr\'{e}sent\'{e}es dans la partie \ref{FCDPartieConnexite}. Le
th\'{e}or\`{e}me \ref{FCD2008FCDSVDgeneralisee} est, dans un certain sens,
optimal, comme le montre le c\'{e}l\`{e}bre exemple suivant~: si l'ensemble
sur lequel $f$ n'est pas d\'{e}rivable n'est pas d\'{e}nombrable, la
conclusion du th\'{e}or\`{e}me est mise en d\'{e}faut.

\begin{exem}
\textbf{L'escalier du diable} (\cite{Rudin1}, voir aussi~:
http://www.mathcurve.com/fractals/ escalierdudiable/escalierdudiable.shtml,
site consult\'{e} le 28 avril 2008).\ D\'{e}finissons $K_{n}$ comme \'{e}tant
la r\'{e}union des $2n$ intervalles ferm\'{e}s $\left[  a/3^{n},\left(
a+1\right)  /3^{n}\right]  $ o\`{u} $a$ est un entier naturel strictement
inf\'{e}rieur \`{a} $3^{n}$ dont l'\'{e}criture en base $3$ comporte au
maximum $n$ chiffres \'{e}gaux \`{a} 0 ou 2. Posons $K=\cap_{n\geq0}K_{n}%
$.\ On rappelle que $K$ est l'ensemble triadique de Cantor.\ Il est compact,
non d\'{e}nombrable, totalement discontinu, de mesure nulle.\ Posons
$f_{0}(x)=x$.\ On d\'{e}finit la fonction $f_{n}$ comme \'{e}tant la fonction
continue, constante sur chacun des intervalles constituant le
compl\'{e}mentaire de l'ensemble $K_{n}$, et affine de pente $(3/2)^{n}$ sur
chacun des intervalles de $K_{n}$. On v\'{e}rifie que, pour tout $n\geq0$
et\ tout $x\in\left[  0,1\right]  $, $\left\vert f_{n+1}(x)-f_{n}%
(x)\right\vert \leq2^{-n}$.\ Ainsi, la s\'{e}rie de fonctions de terme
g\'{e}n\'{e}ral $f_{n+1}-f_{n}$ converge normalement et la suite $\left(
f_{n}\right)  _{n}$ uniform\'{e}ment.\ La limite $f$ est une fonction continue
sur $\left[  0,1\right]  $ croissante comme limite de fonctions croissantes,
v\'{e}rifiant $f(0)=0$, $f(1)=1$.\ Or, $f$ est d\'{e}rivable sur le
compl\'{e}mentaire de l'ensemble non d\'{e}nombrable $K$, et de
d\'{e}riv\'{e}e nulle sur cet ensemble.
\end{exem}

\noindent Nous ne poursuivrons pas ici plus loin ces consid\'{e}rations.\ Ceci
nous entrainerait vers la th\'{e}orie de la mesure qui d\'{e}passe le cadre de
la pr\'{e}sente \'{e}tude. Notons cependant la proximit\'{e} de l'\'{e}tude
des g\'{e}n\'{e}ralisations de la caract\'{e}risation FCD et des conditions
sous lesquelles la relation fondamentale du calcul int\'{e}gral
\[
f(b)-f(a)=\int_{a}^{b}f^{\prime}(t)\,\mathrm{d}t
\]
est vraie. (On peut renvoyer au trait\'{e} d'analyse de W.\ Rudin
\cite{Rudin1} sur ce point.) De la m\^{e}me fa\c{c}on nous n'aborderons pas
les \emph{cas de constance} pour une fonction de la variable complexe holomorphe.

\section{Commentaires et discussion}

\subsection{Un horizon math\'{e}matique ult\'{e}rieur~: op\'{e}rateur de
d\'{e}rivation et caract\'{e}risation FCD}

L'application $D$ qui \`{a} une fonction associe sa d\'{e}riv\'{e}e est un des
premiers exemples d'application lin\'{e}aire issue d'un cadre non
g\'{e}om\'{e}trique que rencontre l'\'{e}l\`{e}ve ou l'\'{e}tudiant.\ Elle
fournit, de plus, un exemple naturel d'application lin\'{e}aire pouvant
\^{e}tre consid\'{e}r\'{e}e sur un espace vectoriel de dimension
infinie.\ Elle est \'{e}galement l'\'{e}l\'{e}ment premier de la th\'{e}orie
des \'{e}quations diff\'{e}rentielles, puisque la caract\'{e}risation FCD est
la plus simple d'entre elles.\ De ce point de vue, cette caract\'{e}risation
est un peu le \textit{juge de paix, }le test de coh\'{e}rence des extensions
de la notion de d\'{e}riv\'{e}e~: on attend d'une d\'{e}riv\'{e}e $D$,
\'{e}tendant la d\'{e}riv\'{e}e usuelle \`{a} un espace $E$ contenant celui
des fonctions usuellement d\'{e}rivables, qu'elle v\'{e}rifie
$D(f)=0\Rightarrow D=$\textrm{ constante}.\smallskip

\noindent Rendons cela plus pr\'{e}cis en nous appuyant sur quelques
exemples.\ Soit $E$ un espace vectoriel ou une alg\`{e}bre munie d'une
d\'{e}riv\'{e}e $D$ interne ($D(E)\subset E$) c'est-\`{a}-dire d'une
op\'{e}ration lin\'{e}aire\ ayant les propri\'{e}t\'{e}s de la d\'{e}riv\'{e}e
usuelle des fonctions. Montrer la caract\'{e}risation FCD revient \`{a}
r\'{e}soudre l'\'{e}quation $D(f)=0$, c'est-\`{a}-dire \`{a} chercher le noyau
$\ker D$ de l'application lin\'{e}aire $D$. Dans le m\^{e}me ordre
d'id\'{e}es, on peut introduire le probl\`{e}me de l'existence de primitives
pour l'op\'{e}rateur $D$ c'est-\`{a}-dire, \'{e}tant donn\'{e} $g\in E$,
d\'{e}terminer s'il existe $f\in E$ tel que $D(f)=g$.\ C'est donc chercher
l'image $\operatorname*{im}D$ de l'application $D$.

\begin{exem}
\label{FCDOpDerExemple1}Soit $E=\mathbb{R}_{n}\left[  x\right]  $, l'ensemble
des fonctions polyn\^{o}mes de degr\'{e} inf\'{e}rieur ou \'{e}gal \`{a} $n$
($n\geq1$). L'alg\`{e}bre $E$ est munie de la d\'{e}riv\'{e}e usuelle des
fonctions que l'on notera $D$. On sait que $\ker D$ est constitu\'{e}e des
fonctions polynomiales constantes (caract\'{e}risation FCD pour les fonctions
polyn\^{o}mes~: en admettant qu'on sache d\'{e}river une fonction
polyn\^{o}me, cette caract\'{e}risation poss\`{e}de alors une
d\'{e}monstration purement alg\'{e}brique.) De la relation $\dim\ker
D+\dim\operatorname*{im}D=n+1$, on sait que $\dim\operatorname*{im}%
D=n$.\ L'op\'{e}rateur $D$ n'est pas surjectif.\ Il est trivial, en effet, de
v\'{e}rifier que la fonction polyn\^{o}me $p_{n}:x\mapsto x^{n}$ n'a pas de
primitives (au sens usuel) dans $\mathbb{R}_{n}\left[  x\right]  $.\ 
\end{exem}

\noindent En revanche, si l'on remplace $\mathbb{R}_{n}\left[  x\right]  $ par
$\mathbb{R}\left[  x\right]  $, l'ensemble des fonctions polyn\^{o}mes,
l'op\'{e}rateur $D$ est surjectif.\ C'est aussi le cas en prenant pour $E$
l'ensemble $\mathrm{C}^{\infty}\left(  \mathbb{R}\right)  $ des fonctions
d\'{e}finies et de classe $\mathrm{C}^{\infty}$ sur le corps des
r\'{e}els\footnote{Voici, d'ailleurs, de bons exemples pour illustrer une des
diff\'{e}rences entre espaces de vectoriels de dimension finie ou infinie.}.
L'analyse moderne a multipli\'{e} les exemples o\`{u} le r\^{o}le de test de
la FCD a \'{e}t\'{e} important.

\begin{exem}
\label{FCDOpDerExemple2}L'espace vectoriel $\mathcal{D}^{\prime}\left(
\mathbb{R}\right)  $ des distributions de Schwartz constitue un autre exemple
d'espace dans lequel tous les objets sont d\'{e}rivables.\ Cet espace contient
l'espace $\mathrm{C}^{0}(\mathbb{R})$ des fonctions continues.\ (En ce sens on
rend d\'{e}rivable les fonctions continues.)\ La d\'{e}riv\'{e}e sur
$\mathcal{D}^{\prime}\left(  \mathbb{R}\right)  $ (not\'{e}e ici $\frac{d}%
{dx}$) prolonge la d\'{e}riv\'{e}e usuelle, dans le sens suivant~: si
$f\in\mathrm{C}^{0}(\mathbb{R})$ est d\'{e}rivable, on a $\frac{d}%
{dx}f=f^{\prime}$. Dans un manuel classique, on note que l'auteur lie bien la
caract\'{e}risation FCD et le probl\`{e}me des primitives, dans une
formulation tr\`{e}s br\`{e}ve. (Il s'agit du b. du th\'{e}or\`{e}me reproduit
dans la figure ci-dessous.)
\end{exem}

%

\begin{figure}
[ptb]
\begin{center}
\includegraphics[
height=2.6602in,
width=5.0298in
]%
{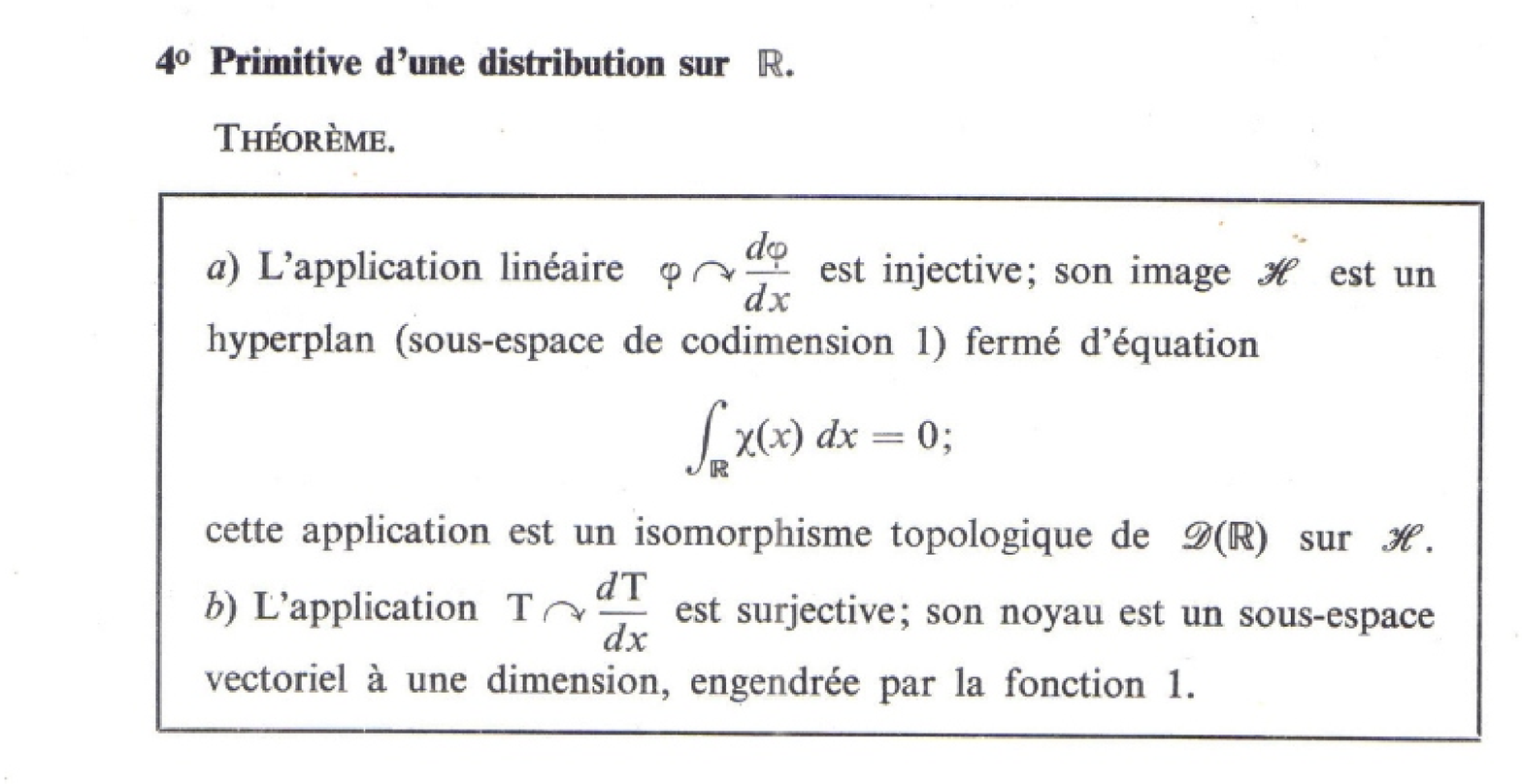}%
\caption{La caract\'{e}risation FCD et l'existence de primitive pour les
distributions d'une variable (Vo-Khac Khoan, 1972)}%
\end{center}
\end{figure}

\section{Un site math\'{e}matique pour la caract\'{e}risation FCD}

Nous utilisons ici, comme dans \cite{SilvyDelcroix1}, une version
l\'{e}g\'{e}rement simplifi\'{e}e et modifi\'{e}e de la notion de site
math\'{e}matique \cite{DuchetErdogan1}.\ Nous organisons ce site autour de
plusieurs champs d'analyse. Les deux premiers sont form\'{e}s d'objets que
nous subdiviserons en deux cat\'{e}gories~:\smallskip

\noindent--~la premi\`{e}re est constitu\'{e}e d'\emph{objets
pr\'{e}construits ou implicites}, de notions protomath\'{e}matiques ou
paramath\'{e}matiques~; ces objets peuvent relever du vocabulaire, de la
logique, de la th\'{e}orie des ensembles ou bien des codages usuels en
math\'{e}matiques.\ Ils peuvent \'{e}galement relever de m\'{e}thodes (au sens
usuels) de d\'{e}monstration~;\ par exemple, pour cette \'{e}tude, figureront
ici deux des strat\'{e}gies de d\'{e}monstration (par l'absurde, par
contrapos\'{e}e) largement utilis\'{e}es. La notion d'intervalle r\'{e}el,
caract\'{e}ris\'{e}s comme \'{e}tant les connexes ou les convexes de
$\mathbb{R}$, est \'{e}galement ici largement implicite de m\^{e}me que celle
de fonction~; l'ensemble de ces objets forme le \emph{substrat} de la
caract\'{e}risation FCD~;\smallskip

\noindent--~la seconde est form\'{e}e des \emph{objets math\'{e}matiques}
centraux de la question math\'{e}matique \'{e}tudi\'{e}e.\ Nous avons choisi
ici de ne faire figurer dans le site que la caract\'{e}risation FCD, dans le
sens ou elle est le coeur de notre sujet. Mais la notion de d\'{e}riv\'{e}e
d'une fonction de la variable r\'{e}elle aurait pu y figurer.\smallskip

\noindent Les \emph{techniques} viennent ensuite. Elles sont entendues ici au
sens de propri\'{e}t\'{e} math\'{e}matique, th\'{e}or\`{e}me en
g\'{e}n\'{e}ral, justifiant une \'{e}tape de la d\'{e}monstration, ici de la
caract\'{e}risation FCD. Elles sont des m\'{e}thodes routini\`{e}res efficaces
et peuvent varier dans les diff\'{e}rentes approches de la question
math\'{e}matique \'{e}tudi\'{e}e. L'in\'{e}galit\'{e} des accroissements finis
(IAF ou IAF') est, par exemple, une des techniques pouvant \^{e}tre
invoqu\'{e}e pour d\'{e}montrer la caract\'{e}risation FCD. La
caract\'{e}risation SVD en est une autre, ainsi que la propri\'{e}t\'{e} de
majoration des accroissements (MAJA). Ces techniques ont pu figurer \`{a} un
moment ou \`{a} un autre au programme des classes du secondaires.\smallskip

\noindent Nous distinguerons enfin plusieurs \emph{niveaux d'analyse
conceptuelle}.\smallskip

\noindent--~Les concepts 1 ou technologie (au sens de Chevallard,
\cite{Chevallard1}) permettent de justifier directement les techniques. Ils
constituent un premier niveau de th\'{e}or\`{e}mes justificatifs~; ainsi,
l'in\'{e}galit\'{e} des acccroissements finis (consid\'{e}r\'{e}e comme
technique de la caract\'{e}risation FCD) peut \^{e}tre justifi\'{e}e par le
th\'{e}or\`{e}me des accroissements finis (TAF), qui en est une technologie.
Nous avons choisi, en fait, de faire figurer ici l'ensemble des
th\'{e}or\`{e}mes accessibles dans l'horizon math\'{e}matique du TAF~:
th\'{e}or\`{e}me de Rolle, propri\'{e}t\'{e}s des fonctions continues.\ (Des
raffinements de ce niveau d'analyse auraient pu se justifier.)\smallskip

\noindent--~les concepts 2 constituent un deuxi\`{e}me niveau de notions ou de
th\'{e}or\`{e}mes justifiant les concepts 1. Clairement,chacune des
technologies \'{e}voqu\'{e}e ci-dessus fait appel \`{a} l'une des
propri\'{e}t\'{e}s $\left(  BS\right)  $, $(SE)$ et $(BW)$ du corps des
r\'{e}els rappel\'{e}es dans l'introduction.\ Nous les avons donc plac\'{e}es
\`{a} ce niveau. Nous notons que ces concepts sont largement implicites dans
le d\'{e}roulement d'un cours de math\'{e}matiques de fin d'\'{e}tudes
secondaires, et qu'elles le restent \'{e}galement souvent en d\'{e}but
d'\'{e}tudes sup\'{e}rieures. (Qui effectue, par exemple, une construction du
corps des r\'{e}els en L1 ou L2 d'un cycle universitaire d'enseignement
sup\'{e}rieur~?) Ce niveau peut \^{e}tre r\'{e}sum\'{e} comme \'{e}tant celui
des propri\'{e}t\'{e}s fondamentales du corps des r\'{e}els.\smallskip

\noindent--~Les concepts 3 et 4 sont ceux de plus haut niveau.\ Ils
constituent les justifications ultimes des concepts mis en oeuvre dans le
probl\`{e}me \'{e}tudi\'{e}.\ Il peut s'agir, selon le cas, de th\'{e}ories
englobantes (qui \'{e}claireraient le probl\`{e}me math\'{e}matique
\'{e}tudi\'{e}), ou bien, \`{a} l'inverse, \^{e}tre au fondement des
math\'{e}matiques mises en jeu dans ce ce probl\`{e}me. A l'inverse des objets
pr\'{e}cons\-truits ou implicites mentionn\'{e}s au d\'{e}but, ils doivent
\^{e}tre explicit\'{e}s (jusqu'\`{a} un certain point) pour rendre claires les
techniques, technologies et diff\'{e}rents concepts mis en jeu dans ce
site.\ Les notions de connexit\'{e}, de compl\'{e}tude et de compacit\'{e},
bref les notions de base en topologie, sont par l'\'{e}clairage qu'elles
apportent aux propri\'{e}t\'{e}s $(BW)$, $\left(  SE\right)  $ et $\left(
BS\right)  $, les principaux champs relevant du troisi\`{e}me niveau
conceptuel. Mais, le calcul diff\'{e}rentiel, les structures alg\'{e}briques
et diff\'{e}rentielles associ\'{e}es sont le cadre ultime dans lesquels se
d\'{e}veloppe la caract\'{e}risation FCD. Ce niveau conceptuel pourrait
\^{e}tre r\'{e}sum\'{e} comme \'{e}tant celui de l'analyse
fonctionnelle.\smallskip

\noindent En r\'{e}sum\'{e} de l'ensemble de l'analyse, nous proposons le site
suivant de la caract\'{e}risation FCD (nous n'avons pas fait figurer les
objets implicites). Les fl\`{e}ches unidirectionnelles indiquent les liens
d'\textit{inclusion dans la justification}, \textquotedblleft\`{a}
interpr\'{e}ter \`{a} peu pr\`{e}s comme
 pertinent pour\textquotedblright\cite{DuchetErdogan1}%
.\ Pour rendre le site plus lisible, des fl\`{e}ches en pointill\'{e}s et
tirets marquent les inclusions sautant un niveau d'analyse.\ Les fl\`{e}ches
ne sont en g\'{e}n\'{e}ral pas indiqu\'{e}es \`{a} l'int\'{e}rieur du m\^{e}me
niveau d'analyse. Enfin, nous avons proc\'{e}d\'{e} \`{a} quelques
regroupements~: le th\'{e}or\`{e}me des accroissements finis, le
th\'{e}or\`{e}me de Rolle, le th\'{e}or\`{e}me de Darboux sont dans un
m\^{e}me cadre, en raison de leur \'{e}quivalence math\'{e}matique~; la
condition n\'{e}cessaire d'extremum (th\'{e}or\`{e}me \ref{FCD2008CNOrdre1})
et les propri\'{e}t\'{e}s des fonctions continues sur les segments
(th\'{e}or\`{e}me \ref{FCD2008theoremeFctContinue1}), comme \'{e}tant les deux
\'{e}l\'{e}ments justificatifs du th\'{e}or\`{e}me de Rolle~; les
propri\'{e}t\'{e}s $(BW)$, $\left(  SE\right)  $ et $\left(  BS\right)  $ en
raison de leur appartenance au groupe des cinq propri\'{e}t\'{e}s
\'{e}quivalentes du corps des r\'{e}els rappel\'{e}es dans l'introduction.%

\begin{figure}
[ptb]
\begin{center}
\includegraphics[
height=3.4662in,
width=4.3154in
]%
{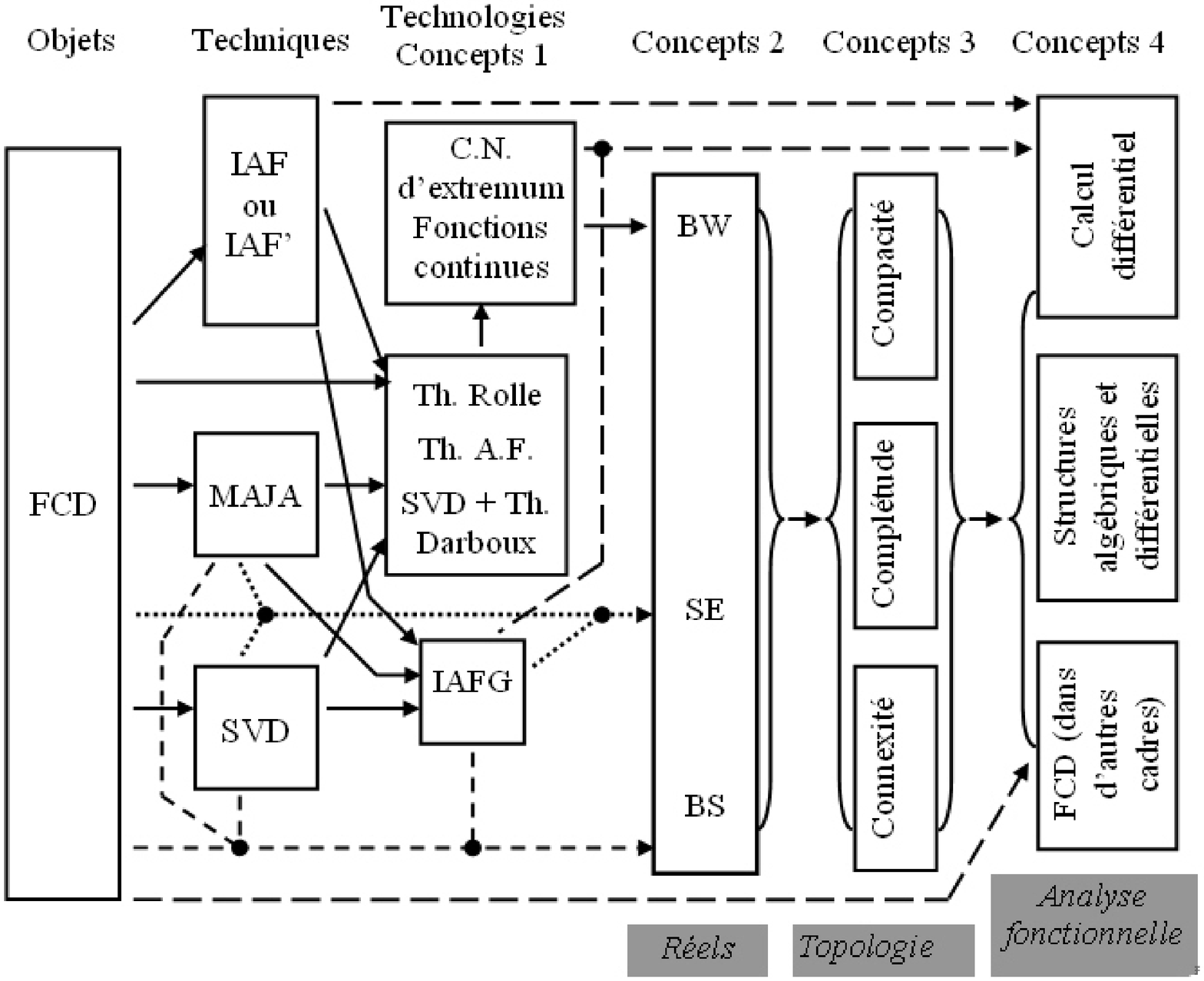}%
\caption{Site de la caract\'{e}risation FCD}%
\end{center}
\end{figure}

\section{Conclusion}

Les consid\'{e}rations math\'{e}matiques et l'analyse en terme de site
pr\'{e}sent\'{e}es plus hauts nous semble montrer que la caract\'{e}risation
des fonctions constantes sur un intervalle par la nullit\'{e} de leur
d\'{e}riv\'{e}e est loin d'\^{e}tre une trivialit\'{e}, tant sur le plan
math\'{e}matique que sur celui de leur enseignement. \smallskip

\noindent D'un point de vue math\'{e}matique, elle plonge ses racines dans les
fondements de l'analyse mo\-derne et s'inscrit dans une classe de
probl\`{e}mes math\'{e}matique plus g\'{e}n\'{e}raux.\ Ses raffinements ont
\'{e}t\'{e} cruciaux dans la bonne compr\'{e}hension de ph\'{e}nom\`{e}nes
li\'{e}s \`{a} la th\'{e}orie de la mesure (l'escalier du diable en est une
des manifestations).\ Elle se trouve, \'{e}galement, assez naturellement
li\'{e}e \`{a} l'expression de conditions optimales pour lesquelles la
relation fondamentale du calcul int\'{e}gral est vraie. Nous mentionnons de
nouveau ici cette propri\'{e}t\'{e} tant son statut semble similaire \`{a}
celui de la caract\'{e}risation FCD~: une \'{e}vidence, mais jusqu'\`{a} quel
point...\smallskip

\noindent Du point de vue de leur enseignement, elle se trouve au coeur d'un
cursus d'analyse math\'{e}matique classique, le mot cursus \'{e}tant pris ici
au sens de culture g\'{e}n\'{e}rale qu'il faut poss\'{e}der en
math\'{e}matique au lyc\'{e}e.\ Ce cursus est situ\'{e} pour ses premiers
horizons \`{a} la fronti\`{e}re entre le cycle terminal et l'enseignement
sup\'{e}rieur mais il se prolonge bien au del\`{a}. Ainsi, le site de la
caract\'{e}risation FCD nous semble-t-il pouvoir \^{e}tre utilis\'{e} tant
dans le cadre de la pr\'{e}paration aux concours du second degr\'{e} que
surtout en formation professionnelle des enseignants pour montrer les
articulations et connexions au sein de ce cursus et expliciter les choix
didactiques effectu\'{e}s par les concepteurs des programmes.\smallskip

\noindent Dans un ordre d'id\'{e}es similaires, cet ensemble de connexions
nous semble pouvoir \'{e}galement \^{e}tre une des raisons explicatives de la
variabilit\'{e} de la situation (pr\'{e}sence, absence,...) des diff\'{e}rents
\'{e}l\'{e}ments de ce site (caract\'{e}risation FCD, SVD, propri\'{e}t\'{e}s
d'accroissements finis,...) dans les programmes du secondaire. (Alors que dans
le m\^{e}me temps le chemin le plus classique pour acc\'{e}der \`{a} la
caract\'{e}risation FCD, par l'in\'{e}galit\'{e} des accroissements finis,
semble faire presque consensus au sein des math\'{e}maticiens.) Mais, par
exemple, le r\^{o}le central de la propri\'{e}t\'{e} de majoration des
accroissements (propri\'{e}t\'{e} MAJA), remis en \'{e}vidence plus haut,
semble sugg\'{e}rer une autre porte d'entr\'{e}e pour l'\'{e}tude des
propri\'{e}t\'{e}s d'accroissements finis en fin de cycle terminal au moins
aussi naturelle que celles adopt\'{e}es dans un pass\'{e}
r\'{e}cent.\smallskip

\noindent Nous esp\'{e}rons que la r\'{e}alisation de ce site math\'{e}matique
de la caract\'{e}risation FCD et les commentaires associ\'{e}s pourront
ajouter un \'{e}l\'{e}ment \`{a} un d\'{e}bat qui nous semble non clos.

\end{document}